\documentclass[review]{elsarticle}

\usepackage{lineno,hyperref}
\usepackage{amssymb,amsmath}
\usepackage[usenames,dvipsnames,svgnames,table]{xcolor}
\modulolinenumbers[250]
\journal{arXiv.org}

\newtheorem{thm}{Theorem}

\newtheorem{prop}[thm]{Proposition}
\newtheorem{corollary}[thm]{Corollary}
\newtheorem{conjecture}[thm]{Conjecture}
\newtheorem{definition}[thm]{Definition}
\newdefinition{rmk}{Remark}
\newproof{pf}{Proof}
\newproof{pot}{Proof of Theorem \ref{thm2}}

\begin{document}

\begin{frontmatter}

\title{Hardy inequality and the construction of infinitesimal operators with non-basis family of eigenvectors}

\author[mymainaddress,mysecondaryaddress]{Grigory M. Sklyar\corref{mycorrespondingauthor}}
\ead{sklar@univ.szczecin.pl}

\author[mythirdaddress]{Vitalii Marchenko}
\cortext[mycorrespondingauthor]{Corresponding author}
\ead{v.marchenko@ilt.kharkov.ua}

\address[mymainaddress]{Institute of Mathematics, University of Szczecin, Wielkopolska 15, 70451, Szczecin, Poland}
\address[mysecondaryaddress]{V.N. Karazin Kharkiv National University, Department of Mathematics and Mechanics, Svobody Sq. 4, 61077, Kharkiv, Ukraine}
\address[mythirdaddress]{Mathematical Division, Institute for Low Temperature Physics and Engineering of NASU, Prospekt Nauky 47, 61103, Kharkiv, Ukraine}

\begin{abstract}
Some special Hilbert spaces are introduced to present the class of infinitesimal operators with complete minimal non-basis family of eigenvectors. The discrete Hardy inequality plays an important role in the proposed approach. The construction complement the results of G.Q.~Xu et al.~\cite{Xu} (2005) and H.~Zwart~\cite{Zwart} (2010) on the Riesz basis property of eigenvectors (eigenspaces) of infinitesimal operators. Our results are extended to the case of operators on some Banach spaces.
\end{abstract}


\end{frontmatter}

\linenumbers

\section{Introduction}
The Hardy inequality (in both its discrete and continuous forms) was discovered at the beginning of XX century and has a lot of applications in various fields of mathematics such as analysis, differential equations, mathematical physics, differential geometry and others~\cite{Kufner1,Kufner2,Opic}. The discrete form of it states~\cite{Hardy1} that if $p>1$ and $\{a_k\}_{k=1}^{\infty}$ is a sequence of nonnegative real numbers, then
\begin{equation}\label{1}
    \sum\limits_{n=1}^{\infty} \left(\frac{1}{n} \sum\limits_{k=1}^{n} a_k\right)^p \leq \left(\frac{p}{p-1}\right)^p \sum\limits_{n=1}^{\infty} a_n^p.
\end{equation}
We notice that G.H.~Hardy came to the discovery of inequality~(\ref{1}) when he tried to obtain an elementary proof of (the weak form of) the Hilbert inequality~\cite{Hardy2,Kufner1}. The latter inequality was discovered when D.~Hilbert studied the solutions to a certain integral equations, see~\cite{Hilbert}, and the weak form of it asserts that if $\{a_n\}_{n=1}^{\infty},\{b_n\}_{n=1}^{\infty}\in \ell_2$ and $a_n\geq 0,\:b_n\geq 0,$ then the double series $\sum\limits_{n=1}^{\infty}\sum\limits_{m=1}^{\infty} \frac{a_m b_n}{m+n}$ converges. We must also stress an essential contribution of such mathematicians as M.~Riesz, E.~Landau and I.~Schur to the development of~(\ref{1}). For more details see, e.g.,~\cite{Kufner1}.

In the following we denote by $H$ a separable Hilbert space with norm $\|\cdot\|$ and scalar product $\langle\cdot,\cdot\rangle$. The central concept of this paper is the concept of $C_0$-group of linear bounded operators. If $A$ is an infinitesimal generator of the $C_0$-group on $H$, then the abstract Cauchy problem
\begin{equation}\label{Cauchy problem}
    \left \{
\begin{array}{l}
\dot{x}(t)=Ax(t),\quad t \in\mathbb{R},\\
x(0)=x_0,\\
\end{array}
\right .
\end{equation}
on $H$ is well-posed in the sense of J.~van~Neerven (for each $x_0\in D(A)$ the problem~(\ref{Cauchy problem}) has a unique classical solution) and $\rho(A)\neq \emptyset$, see~\cite{vanNeerven}.
A remarkable result in the spectral theory of $C_0$-semigroups in Hilbert spaces was obtained in~\cite{Xu,Zwart}, where the Riesz basis property for eigenvectors of certain class of generators of $C_0$-groups was established. We give here a formulation of the main result of these works in some particular case.
\begin{thm}[\cite{Zwart}]\label{Remarkable result}
Let $A$ be the generator of the $C_0$-group on $H$ with eigenvalues $\{\lambda_n\}_{n=1}^{\infty}$ (counting with multiplicity) and the corresponding (normalized) eigenvectors $\{e_n\}_{n=1}^{\infty}$. If the following two conditions hold,
\begin{enumerate}
\item $\overline{Lin}\{e_n\}_{n=1}^{\infty}=H;$
\item The point spectrum has a uniform gap, i.e.,
\begin{equation}\label{2}
   \inf\limits_{n\neq m} |\lambda_n-\lambda_m|>0,
\end{equation}
\end{enumerate}
then $\{e_n\}_{n=1}^{\infty}$ forms a Riesz basis of $H$.
\end{thm}

We note that Theorem~\ref{Remarkable result} follows from the main result of~\cite{Xu}, but the approach used by H.~Zwart in~\cite{Zwart} essentially differs from the one proposed by G.Q.~Xu and S.P.~Yung in~\cite{Xu}. The proof of Theorem~\ref{Remarkable result} in~\cite{Zwart} is based on Carleson's interpolation theorem~\cite{Garnett} and on the fact that each generator of a $C_0$-group on $H$ has a bounded $\mathcal{H}^{\infty}$-calculus on a strip~\cite{Boyadzhiev,Haase1,Haase2}.

Let the eigenvalues $\{\lambda_n\}_{n=1}^{\infty}$ of the generator $A$ of the $C_0$-group on $H$ can be grouped into $K$ sets $\{\lambda_{n,1}\}_{n=1}^{\infty},\{\lambda_{n,2}\}_{n=1}^{\infty},\dots,\{\lambda_{n,K}\}_{n=1}^{\infty}$ with $\inf\limits_{n\neq m} |\lambda_{n,k}-\lambda_{m,k}|>0,$ $k=1,\dots,K,$ and the span of the generalized eigenvectors of $A$ is dense. Then, as it is shown in~\cite{Zwart}, there exists a sequence of spectral projections $\{P_n\}_{n=1}^{\infty}$ of $A$ such that $\{P_n H\}_{n=1}^{\infty}$ forms a Riesz basis of subspaces in $H$ with $\max\limits_{n} \dim P_n H\leq K.$

The main goal of our work is to show that the assumption~(\ref{2}) in Theorem~\ref{Remarkable result} is essential, i.e. if we omit this assumption or even weaken it, the statement of Theorem~\ref{Remarkable result} becomes false. The case when $\{\lambda_n\}_{n=1}^{\infty}$ can be decomposed into $K$ sets, with every set satisfying~(\ref{2}), was considered in~\cite{Zwart}. So we consider the case when $\{\lambda_n\}_{n=1}^{\infty}$ does not satisfy~(\ref{2}) and, moreover, cannot be decomposed into $K$ sets, with every set satisfying~(\ref{2}). We present the construction of the generator $A$ of the $C_0$-group with eigenvalues $\{\lambda_n\}_{n=1}^{\infty}$ and complete minimal non-basis family of eigenvectors (Theorem~\ref{not-basis-ln}). It must be emphasized that the application of the discrete Hardy inequality~(\ref{1}) (for $p=2$) is crucial in our approach.

Furthermore, we give some generalizations of our result and present a class of infinitesimal operators with complete minimal non-basis family of eigenvectors (Theorem~\ref{not-basis_fk}). For this purpose we consider special classes of real sequences $\mathcal{S}_k$, $k\in\mathbb{N}$, present classes of Hilbert spaces $H_k\left(\{e_n\}\right)$, $k\in\mathbb{N}$, depending on $H$ and on a chosen Riesz basis $\{e_n\}_{n=1}^{\infty}$ of $H$, and prove that $\{e_n\}_{n=1}^{\infty}$ do not form a Schauder basis of $H_k\left(\{e_n\}\right)$, $k\in\mathbb{N}$.

Further on, we propose certain development of our approach. Namely, we use the same idea for the construction of infinitesimal operators, acting in certain Banach spaces, with complete minimal non-basis family of eigenvectors. This construction is essentially based on Hardy inequality~(\ref{1}) for $p>1$.
In Theorem~\ref{not-basis-ell_p_fk} we present a class of infinitesimal operators on Banach spaces with complete minimal non-basis family of eigenvectors. To this end we use classes of real sequences $\mathcal{S}_k$, $k\in\mathbb{N}$, introduce certain classes of Banach sequence spaces $\ell_{p,k}\left(\{e_n\}\right),$ $p\geq 1,$ $k\in\mathbb{N}$, depending on given $\ell_p$ space and on arbitrary chosen symmetric basis $\{e_n\}_{n=1}^{\infty}$ of $\ell_p$, and show that $\{e_n\}_{n=1}^{\infty}$ do not form a Schauder basis of $\ell_{p,k}\left(\{e_n\}\right).$ The concept of symmetric basis was first introduced and studied by I.~Singer~\cite{Singer2} in connection with S.~Banach's hyperplane problem. For various properties of symmetric bases see, e.g.,~\cite{Lindenstrauss,Singer}.
\section{Auxiliary constructions and preliminary results}
\subsection{Spaces $H_k\left(\{e_n\}\right)$, $k\in\mathbb{N}$}
We begin with the following definition.
\begin{definition}\label{right shift}
Let $E$ be a Banach space with basis $\{e_n\}_{n=1}^{\infty}$. Then the operator $T$ defined on $E$ by $T e_n=e_{n+1},\: n\in\mathbb{N},$ will be called the right shift operator associated with the basis  $\{e_n\}_{n=1}^{\infty}$.
\end{definition}
Suppose that $\{e_n\}_{n=1}^{\infty}$ is an arbitrary Riesz basis of $H$ and $T$ is the right shift operator associated with $\{e_n\}_{n=1}^{\infty}$. We introduce the following spaces,
$$H_k^0\left(\{e_n\}\right)=\left\{ x\in H:\: \|x\|_k =\left\|\left(I-T\right)^k x\right\|\right\},\: k\in\mathbb{N}.$$
We notice that $0\in \sigma\left((I-T)^k\right)$ for any $k$. Hence, $H_k^0\left(\{e_n\}\right)$ is a normed linear space, but not complete.
By $H_k\left(\{e_n\}\right)$ we denote the completion of $H_k^0\left(\{e_n\}\right)$ in the norm $\|\cdot\|_k$. Using the characteristic property of Riesz basis we obtain
$$H=\left\{x=\sum\limits_{n=1}^{\infty}c_n e_n:\: \{c_n\}_{n=1}^{\infty}\in \ell_2\right\}.$$
Next, we observe that
\begin{align*}
\|x\|_k &=\left\|\left(I-T\right)^k \sum\limits_{n=1}^{\infty}c_n e_n\right\|\\
&=\left\|\sum\limits_{n=1}^{\infty}c_n \left(e_{n} -C_{k}^{1} e_{n+1}+\dots + (-1)^{k-1} C_{k}^{k-1} e_{n+k-1} +(-1)^{k} e_{n+k}\right)\right\|\\
&=\left\|\sum\limits_{n=1}^{\infty}\left(c_{n} -C_{k}^{1} c_{n-1}+\dots+  (-1)^{k+1} C_{k}^{k-1} c_{n-k+1} +(-1)^{k} c_{n-k}\right) e_n\right\|,
\end{align*}
where for $j\in\mathbb{N}$ we set $c_{1-j}=0$.
The last norm is finite if and only if the condition
$$\sum\limits_{n=1}^{\infty} \left|c_{n} -C_{k}^{1} c_{n-1}+\dots +(-1)^{k} c_{n-k}\right|^2<\infty
$$
holds. Consequently, for each $k\in\mathbb{N},$ the space $H_k\left(\{e_n\}\right)$ consists of formal series $x=(\mathfrak{f})\sum\limits_{n=1}^{\infty}c_n e_n$ with the property
$$\left\{c_{n} -C_{k}^{1} c_{n-1}+\dots +(-1)^{k} c_{n-k}\right\}_{n=1}^{\infty}\in \ell_2.$$

It turns out that $H_k\left(\{e_n\}\right)$ is a Hilbert space with a norm
\begin{equation}\label{norm}
  \|x\|_k=\left\|(\mathfrak{f})\sum\limits_{n=1}^{\infty}c_n e_n\right\|_k=\left\|\sum\limits_{n=1}^{\infty}\left(c_{n} -C_{k}^{1} c_{n-1}+\dots +(-1)^{k} c_{n-k}\right) e_n\right\|,
\end{equation}
$x\in H_k\left(\{e_n\}\right),$
and a scalar product
$$\langle x,y \rangle_{k}=\left\langle \left(I-T\right)^k x,\left(I-T\right)^k y \right\rangle,\: x,y\in H_k\left(\{e_n\}\right).$$

For example, given any $\alpha\in [0,\frac{1}{2})$, we have $(\mathfrak{f})\sum\limits_{n=1}^{\infty} n^{\alpha} e_n\in H_1\left(\{e_n\}\right)$. Indeed, for $\alpha=0$ this fact is obvious. For $\alpha\in (0,\frac{1}{2})$ we have that
$$n^{\alpha}-(n-1)^{\alpha}\sim c_{\alpha} n^{\alpha-1}$$
as $n\rightarrow \infty.$
Hence, $\left\{n^{\alpha}-(n-1)^{\alpha}\right\}_{n=1}^{\infty}\in \ell_2.$

Concerning the inner product $\langle\cdot,\cdot\rangle_k$ we can say more.
If $\{e_n\}_{n=1}^{\infty}$ is an orthonormal basis of $H$, $x=(\mathfrak{f})\sum\limits_{n=1}^{\infty}c_n e_n$, $y=(\mathfrak{f})\sum\limits_{n=1}^{\infty}d_n e_n\in H_k\left(\{e_n\}\right)$, then
$$\langle x,y \rangle_{k}= \sum\limits_{n=1}^{\infty} (c_{n} -C_{k}^{1} c_{n-1}+\dots +(-1)^{k} c_{n-k})(\overline{d}_{n} -C_{k}^{1} \overline{d}_{n-1}+\dots +(-1)^{k} \overline{d}_{n-k}).$$

We note that, in particular case of $H=\ell_2$ and when $\{e_n\}_{n=1}^{\infty}$ is the canonical basis of $\ell_2,$ $H_k\left(\{e_n\}\right)=\ell_2(\Delta^k).$
A sequence space $\ell_2(\Delta^k)$ is the space consisting of all sequences whose $k^{th}$ order differences are $2$-absolutely summable, with norm $\|x\|_{\ell_2(\Delta^k)}=\|\Delta^k x\|_{\ell_2},$ where $\Delta$ is a difference operator, i.e.
$$\Delta=\left( \begin{array}{ccccc} 1 & 0 & 0 & 0 & \dots\\ -1 & 1 & 0 & 0 & \dots\\ 0 & -1 & 1 & 0 & \dots\\ 0 & 0 & -1 & 1 & \dots\\ \vdots & \vdots & \vdots & \vdots & \ddots\end{array} \right),$$
see~\cite{Altay,Imaninezhad} and earlier paper~\cite{Basar}, where only the case $k=1$ is considered. In other words, $\ell_2(\Delta^k)=\left\{x=\{\alpha_n\}_{n=1}^{\infty}:\:\Delta^k x \in \ell_2\right\}.$ Hence,
\begin{equation}\label{H_k}
  H_k\left(\{e_n\}\right)=\left\{x=(\mathfrak{f})\sum\limits_{n=1}^{\infty}c_n e_n:\: \{c_n\}_{n=1}^{\infty}\in \ell_2(\Delta^k)\right\},\: k\in\mathbb{N}.
\end{equation}
Thus, our class of spaces $H_k\left(\{e_n\}\right)$ is analogous to $\ell_2(\Delta^k)$, studied in~\cite{Altay,Basar,Imaninezhad}. Moreover, we note that $\ell_2(\Delta^k)$ naturally arises as a completion of $\left(\ell_2\right)_k^0\left(\{e_n\}\right)$, where $\{e_n\}_{n=1}^{\infty}$ is a canonical basis of $\ell_2.$
The following proposition includes some properties of the space $H_k\left(\{e_n\}\right)$.
\begin{prop}\label{prop}
The space $H_k\left(\{e_n\}\right)$ has the following properties.
\begin{enumerate}
\item $\overline{Lin}\{e_n\}_{n=1}^{\infty}=H_k\left(\{e_n\}\right);$
\item $\{e_n\}_{n=1}^{\infty}$ does not form a basis of $H_k\left(\{e_n\}\right);$
\item $\{e_n\}_{n=1}^{\infty}$ has a unique biorthogonal system $$\left\{\chi_n=\left(I-T\right)^{-k} \left(I-T^{\ast}\right)^{-k} e_n^{\ast}\right\}_{n=1}^{\infty}$$
      in $H_k\left(\{e_n\}\right)$, where $\langle e_{n},e_{m}^{\ast} \rangle = \delta_{n}^{m}$;
\item $\{\chi_n\}_{n=1}^{\infty}$ is uniformly minimal sequence in $H_k\left(\{e_n\}\right)$ while the sequence $\{e_n\}_{n=1}^{\infty}$ is minimal but not uniformly minimal in
      $H_k\left(\{e_n\}\right);$
\item $H\subset H_1\left(\{e_n\}\right)\subset H_2\left(\{e_n\}\right)\subset H_3\left(\{e_n\}\right)\subset \dots$ ;
\item $H_k\left(\{e_n\}\right)$ is a separable Hilbert space, isomorphic to $\ell_2$;
\item $H_k\left(\{e_n\}\right)$ has an orthonormal basis;
\item $L=\left\{x=(\mathfrak{f})\sum\limits_{n=1}^{\infty}c_n e_n\in H_k\left(\{e_n\}\right):\: \{c_n\}_{n=1}^{\infty}\in \ell_2(\Delta^k)\cap c_0\right\}$, where $c_0$ is the space of sequences $\{\alpha_n\}_{n=1}^{\infty}$ with $\lim\limits_{n\rightarrow \infty} \alpha_n=0$, is not a subspace of $H_k\left(\{e_n\}\right)$.
\end{enumerate}
\end{prop}

\begin{pf} \textit{1.} Follows from the fact that only zero is orthogonal to all the $e_n$'s with respect to the scalar product $\langle\cdot,\cdot\rangle_k$.

\textit{2.} Is a consequence of the fact that  $(\mathfrak{f})\sum\limits_{n=1}^{\infty} e_n\in H_k\left(\{e_n\}\right)$ can not be represented by the series $\sum\limits_{n=1}^{\infty} c_n e_n$, convergent in $H_k\left(\{e_n\}\right)$, since $\inf\limits_{n} \|e_n\|_k >0.$

\textit{3.} Clearly, $\langle e_{n},\chi_j \rangle_k = \left\langle \left(I-T\right)^{k} e_{n}, \left(I-T\right)^{k} \left(I-T\right)^{-k} \left(I-T^{\ast}\right)^{-k} e_j^{\ast} \right\rangle =\delta_{n}^{j}$, and the uniqueness of $\{\chi_n\}_{n=1}^{\infty}$ follows from \textit{1.}

\textit{4.} Is true since $\sup\limits_{n} \|e_n\|_k <\infty$ while $\sup\limits_{n} \|\chi_n\|_k =\infty$ \cite{Keown}.

\textit{5.} Follows from the chain of inclusions $\ell_2\subset \ell_2(\Delta)\subset \ell_2(\Delta^2)\subset \ell_2(\Delta^3)\subset\dots$~\cite{Imaninezhad}.

\textit{6.} We fix $x=(\mathfrak{f})\sum\limits_{n=1}^{\infty}c_n e_n\in H_k\left(\{e_n\}\right)$ and denote $(\Delta^k c)_n=c_{n} -C_{k}^{1} c_{n-1}+\dots +(-1)^{k} c_{n-k}.$
Combining (\ref{norm}) with the property of Riesz basis in $H$, we obtain the following inequality,
$$m \sum\limits_{n=1}^{\infty}\left|(\Delta^k c)_n\right|^2 \leq\|x\|_k^2\leq M \sum\limits_{n=1}^{\infty}\left|(\Delta^k c)_n\right|^2,$$
which generates an isomorphism between $H_k\left(\{e_n\}\right)$ and $\ell_2(\Delta^k)$. And, since $\ell_2(\Delta^k)$ is isometrically isomorphic to $\ell_2$ \cite{Imaninezhad}, $H_k\left(\{e_n\}\right)$ is isomorphic to $\ell_2.$ Hence, $H_k\left(\{e_n\}\right)$ is a separable space.

\textit{7.} Is a consequence of well-known fact that every separable Hilbert space has an orthonormal basis.

\textit{8.} The proof is based on the fact that $\ell_2(\Delta^k)\cap c_0$ is not closed in $\ell_2(\Delta^k)$.
\hfill$\Box$
\end{pf}

For example, if $\{e_n\}_{n=1}^{\infty}$ is an orthonormal basis in $H$, then it is clear that the sequence
$$\left\{\left(I-T\right)^{-k} e_n\right\}_{n=1}^{\infty}$$
forms an orthonormal basis of $H_k\left(\{e_n\}\right)$.
From the other hand, it is interesting to construct a (bounded) non-Riesz basis of $H_k\left(\{e_n\}\right)$.
We recall that the first example of bounded non-Riesz basis appeared only in~1948 and it was given by K.I.~Babenko in~\cite{Babenko}. He showed that for every $\alpha\in\left(-\frac{1}{2},0\right) \cup \left(0,\frac{1}{2}\right)$, the system of functions $\{|t|^{\alpha} e^{int}\}_{n=-\infty}^{\infty}$ forms a bounded non-Riesz basis in $L_2 (-\pi,\pi)$. This example was later generalized by V.F.~Gaposhkin~\cite{Gaposhkin} and operators generating non-Riesz bases in $H$ were studied by A.M.~Olevskii in~\cite{Olevskii}.

\subsection{Spaces $\ell_{p,k}\left(\{e_n\}\right),$ $p\geq 1,$ $k\in\mathbb{N}$}

Similarly to the above we introduce the space $\ell_{p,k}\left(\{e_n\}\right)$ as a completion of the space
$$\ell_{p,k}^0\left(\{e_n\}\right)=\left\{ x\in \ell_p:\: \|x\|_k =\left\|\left(I-T\right)^k x\right\|\right\},\: k\in\mathbb{N},$$
where $\{e_n\}_{n=1}^{\infty}$ is a symmetric basis of $\ell_p,$ $p\geq 1,$ and $T$ is the right shift operator associated with $\{e_n\}_{n=1}^{\infty}$. It is known that spaces $\ell_p$, $p\geq 1$, have a unique, up to equivalence, symmetric basis~\cite{Lindenstrauss} and it is equivalent to the canonical one. Thus we arrive at the following assertion.
\begin{prop}\label{sym}
Let $\{e_n\}_{n=1}^{\infty}$ be a symmetric basis of $\ell_p,\:p\geq 1$. Then there exist constants $M\geq m >0$ such that for each $x=\sum\limits_{n=1}^{\infty} c_n e_n \in \ell_p,$
$$m\sum\limits_{n=1}^{\infty} |c_n|^p \leq \|x\|^p \leq M \sum\limits_{n=1}^{\infty} |c_n|^p.$$
\end{prop}
Proposition~\ref{sym} says that the class of Riesz bases in $\ell_2$ coincides with the class of symmetric bases.
Using Proposition \ref{sym} we obtain that
$$\ell_p=\left\{x=\sum\limits_{n=1}^{\infty}c_n e_n:\: \{c_n\}_{n=1}^{\infty}\in \ell_p\right\}.$$
Consequently, for each $k\in\mathbb{N},$
$$\ell_{p,k}\left(\{e_n\}\right)=\left\{x=(\mathfrak{f})\sum\limits_{n=1}^{\infty}c_n e_n:\: \left\{c_{n} -C_{k}^{1} c_{n-1}+\dots +(-1)^{k} c_{n-k}\right\}_{n=1}^{\infty}\in \ell_p\right\}$$
is a Banach sequence space with norm
$$\|x\|_k=\left\|(\mathfrak{f})\sum\limits_{n=1}^{\infty}c_n e_n\right\|_k=\left\|\sum\limits_{n=1}^{\infty}\left(c_{n} -C_{k}^{1} c_{n-1}+\dots+(-1)^{k} c_{n-k}\right) e_n\right\|,\: x\in \ell_{p,k}\left(\{e_n\}\right).$$

Except for the case $p=2$, the space $\ell_{p,k}\left(\{e_n\}\right),$ $k\in\mathbb{N},$ is not an inner product space and, hence, not a Hilbert space. On the other hand, $\ell_{2,k}\left(\{e_n\}\right),$ $k\in\mathbb{N},$ is a Hilbert space since $\ell_{2,k}\left(\{e_n\}\right)=H_k\left(\{e_n\}\right),$ where $H=\ell_2$ and $\{e_n\}_{n=1}^{\infty}$ is a Riesz basis of $\ell_2.$
We also note that, if $\{e_n\}_{n=1}^{\infty}$ denotes the canonical basis of $\ell_p,$ then $\ell_{p,k}\left(\{e_n\}\right)=\ell_p(\Delta^k)$, where $\ell_p(\Delta^k)$ is the space consisting of all sequences whose $k^{th}$ order differences are $p$-absolutely summable, with norm $\|x\|_{\ell_p(\Delta^k)}=\|\Delta^k x\|_{\ell_p},$ see~\cite{Altay,Imaninezhad,Basar} for details. In other words, $\ell_p(\Delta^k)=\left\{x=\{\alpha_n\}_{n=1}^{\infty}:\:\Delta^k x \in \ell_p\right\}.$ It follows that
$$\ell_{p,k}\left(\{e_n\}\right)=\left\{x=(\mathfrak{f})\sum\limits_{n=1}^{\infty}c_n e_n:\: \{c_n\}_{n=1}^{\infty}\in \ell_p(\Delta^k)\right\},\: k\in\mathbb{N}.$$

In the following proposition we collect some properties of the space $\ell_{p,k}\left(\{e_n\}\right).$
\begin{prop}\label{props}
Let $\{e_n\}_{n=1}^{\infty}$ be a symmetric basis of $\ell_p$, $p\geq 1,$ and $k\in\mathbb{N}$. Then the following statements are true.
\begin{enumerate}
\item If $p>1$, then $\overline{Lin}\{e_n\}_{n=1}^{\infty}=\ell_{p,k}\left(\{e_n\}\right);$
\item $\{e_n\}_{n=1}^{\infty}$ does not form a basis of $\ell_{p,k}\left(\{e_n\}\right);$
\item If $p>1$, then $\{e_n\}_{n=1}^{\infty}$ has a unique biorthogonal system $$\left\{\chi_n=\left(I-T\right)^{-k} \left(I-T^{\ast}\right)^{-k} e_n^{\ast}\right\}_{n=1}^{\infty}$$
      in $\left(\ell_{p,k}\left(\{e_n\}\right)\right)^{\ast}$, where $\{e_n^{\ast}\}_{n=1}^{\infty}$ is biorthogonal to $\{e_n\}_{n=1}^{\infty}$ basis of $\ell_q$, where $\frac{1}{p}+\frac{1}{q}=1$;
\item If $p>1$, then $\{\chi_n\}_{n=1}^{\infty}$ is uniformly minimal sequence in $\left(\ell_{p,k}\left(\{e_n\}\right)\right)^{\ast}$ while the sequence $\{e_n\}_{n=1}^{\infty}$ is minimal but not uniformly minimal in $\ell_{p,k}\left(\{e_n\}\right);$
\item $\ell_p\subset \ell_{p,1}\left(\{e_n\}\right)\subset \ell_{p,2}\left(\{e_n\}\right)\subset \ell_{p,3}\left(\{e_n\}\right)\subset \dots$ ;
\item $\ell_{p,k}\left(\{e_n\}\right)$ is a separable Banach sequence space, isomorphic to $\ell_p$;
\item $L=\left\{x=(\mathfrak{f})\sum\limits_{n=1}^{\infty}c_n e_n\in \ell_{p,k}\left(\{e_n\}\right):\: \{c_n\}_{n=1}^{\infty}\in \ell_p(\Delta^k)\cap c_0\right\}$ is not a subspace of $\ell_{p,k}\left(\{e_n\}\right)$.
\end{enumerate}
\end{prop}
The proof of Proposition~\ref{props} is similar to the proof of Proposition~\ref{prop}.
\section{The construction of infinitesimal operators with non-basis family of eigenvectors on Hilbert spaces}
\subsection{Infinitesimal operators on $H_1\left(\{e_n\}\right)$}
In the following, by $[X]$ we denote the space of all bounded linear operators on a Banach space~$X.$ We define the operator $A: H_1\left(\{e_n\}\right) \supset D(A) \mapsto H_1\left(\{e_n\}\right)$ by formula
\begin{equation}\label{A}
  A x=A (\mathfrak{f})\sum\limits_{n=1}^{\infty} c_n e_n= (\mathfrak{f})\sum\limits_{n=1}^{\infty} \lambda_n c_n e_n,
\end{equation}
where $\{\lambda_n\}_{n=1}^{\infty}$ is an unbounded sequence on the complex plane, and domain
\begin{equation}\label{Domain_A}
    D(A)=\left\{x= (\mathfrak{f})\sum\limits_{n=1}^{\infty} c_n e_n \in H_1\left(\{e_n\}\right):\:\{\lambda_n c_n\}_{n=1}^{\infty}\in \ell_2(\Delta)\right\}.
\end{equation}
An example of the generator of unbounded $C_0$-group with eigenvectors that do not form a basis is given by the following theorem.
\begin{thm}\label{not-basis-ln}
The operator $A$ defined by (\ref{A}) with domain (\ref{Domain_A}), where $\lambda_n=i\ln n$, $n\in\mathbb{N}$,
generates the $C_0$-group $\{e^{A t}\}_{t\in \mathbb{R}}$ on $H_1\left(\{e_n\}\right)$, which acts for every $t\in\mathbb{R}$ by the formula
\begin{equation}\label{Solution}
e^{A t} x=e^{A t}(\mathfrak{f})\sum\limits_{n=1}^{\infty} c_n e_n =(\mathfrak{f})\sum\limits_{n=1}^{\infty} e^{i t \ln n} c_{n} e_n.
\end{equation}
\end{thm}
\begin{pf} Recall that $\{e_n\}_{n=1}^{\infty}$ is a Riesz basis of $H$, hence there exist
constants $M\geq m>0$ such that for every $x=\sum\limits_{n=1}^{\infty} \alpha_n \phi_n\in H$ we have
\begin{equation}\label{Riesz}
  m \|x\|^2 \leq \sum\limits_{n=1}^{\infty} |\alpha_n|^2 \leq M \|x\|^2.
\end{equation}

Also we recall that the norm of element $x=(\mathfrak{f})\sum\limits_{n=1}^{\infty}c_n e_n\in H_1\left(\{e_n\}\right)$ is computed as follows,
\begin{equation}\label{1-norm}
  \|x\|_1=\left\|(\mathfrak{f})\sum\limits_{n=1}^{\infty}c_n e_n\right\|_1=\left\|\sum\limits_{n=1}^{\infty}\left(c_{n} - c_{n-1}\right) e_n\right\|,
\end{equation}
where $c_0=0.$ To prove that  $(\mathfrak{f})\sum\limits_{n=1}^{\infty} e^{i t \ln n} c_{n} e_n \in H_1\left(\{e_n\}\right)$ let us
consider $x=(\mathfrak{f})\sum\limits_{n=1}^{\infty} c_n e_n\in H_1\left(\{e_n\}\right)$ and fix $t\in\mathbb{R}$.
Then, by~(\ref{1-norm}) we will have
\begin{align*}
&\left\|(\mathfrak{f})\sum\limits_{n=1}^{\infty} e^{i t \ln n} c_{n} e_n\right\|_{1}=\left\|c_1 e_1 + \sum\limits_{n=2}^{\infty}\left(e^{i t \ln n} c_{n} -e^{i t \ln (n-1)} c_{n-1}\right) e_n\right\|\\
&=\left\|c_1 e_1 + \sum\limits_{n=2}^{\infty} ( e^{i t \ln n} c_{n} - e^{i t \ln n} c_{n-1}+e^{i t \ln n} c_{n-1}-e^{i t \ln (n-1)} c_{n-1}) e_n\right\|\\
&\leq\left\|\sum\limits_{n=1}^{\infty} e^{i t \ln n} (c_n-c_{n-1}) e_n\right\| + \left\|\sum\limits_{n=2}^{\infty} (e^{i t \ln n}-e^{i t \ln (n-1)}) c_{n-1} e_n\right\|=\Xi_1 + \Xi_2(t).
\end{align*}
The first term is estimated immediately by using~(\ref{Riesz}),
\begin{equation}\label{Ksi1}
\Xi_1^2\leq \frac{1}{m} \sum\limits_{n=1}^{\infty} \left| (c_n-c_{n-1}) \right|^2 \leq \frac{M}{m} \left\| \sum\limits_{n=1}^{\infty} (c_n-c_{n-1}) e_n\right\|^2 = \frac{M}{m} \|x\|_1^2.
\end{equation}


Further, since $e^{i t \ln n} - e^{i t \ln (n-1)}=e^{i t \ln n}\left(1-e^{it \ln\left(1-\frac{1}{n}\right)}\right)$ is valid for all $n\geq 2$, then, by using~(\ref{Riesz}) we will have
$$\Xi_2^2(t)\leq \frac{1}{m} \sum\limits_{n=2}^{\infty} \left|\left(1-e^{it \ln\left(1-\frac{1}{n}\right)}\right) c_{n-1} \right|^2 = \frac{1}{m} \sum\limits_{n=2}^{\infty} \frac{n^2 \left|1-e^{it \ln\left(1-\frac{1}{n}\right)} \right|^2}{t^2} \frac{t^2}{n^2} |c_{n-1}|^2.$$
We introduce the notation $\xi_n(t)=\frac{n^2 \left|1-e^{it \ln\left(1-\frac{1}{n}\right)} \right|^2}{t^2}$. Note that, for $r\geq 2$ we have $r \left| \ln\left(1-\frac{1}{r}\right)\right|\leq 2$, and for all  $s\in\mathbb{R}$ the following inequalities are true,
$$\sin^2 s\leq s^2,\qquad \left(1-\cos s \right)^2 \leq s^2.$$
Consequently, for all $t\in\mathbb{R}$ and all $n\geq 2$ we have
\begin{align*}
\xi_n(t) &=\frac{n^2}{t^2}\left( \left(1-\cos\left(t \ln\left(1-\frac{1}{n}\right) \right) \right)^2 + \sin^2\left(t \ln\left(1-\frac{1}{n}\right)\right) \right) \\
&\leq \frac{n^2}{t^2} \left( 2\left( t \ln\left(1-\frac{1}{n}\right) \right)^2\right)=2 n^2 \left(\ln\left(1-\frac{1}{n}\right) \right)^2 \leq 8.
\end{align*}

Thus,
$$
\Xi_2^2(t)\leq \frac{8 t^2}{m} \sum\limits_{n=2}^{\infty} \frac{|c_{n-1}|^2}{n^2}.
$$
Next, to estimate $\sum\limits_{n=2}^{\infty} \frac{|c_{n-1}|^2}{n^2}$ by $C\sum\limits_{n=1}^{\infty} |c_n-c_{n-1}|^2$ we remark that
$$c_{n-1}= \sum\limits_{j=1}^{n-1} \left(c_j-c_{j-1}\right),\quad n\geq 2,$$
and at the beginning we estimate $\Xi_2^2(t)$ as follows,
\begin{align*}
\Xi_2^2(t) &\leq \frac{8 t^2}{m} \sum\limits_{n=2}^{\infty} \frac{|c_{n-1}|^2}{n^2}\leq \frac{8 t^2}{m} \sum\limits_{n=2}^{\infty} \frac{1}{n^2} \left|\sum\limits_{j=1}^{n-1} \left(c_j-c_{j-1}\right)\right|^2\\
& \leq \frac{8 t^2}{m} \sum\limits_{n=2}^{\infty}
\left(\frac{1}{n} \sum\limits_{j=1}^{n-1}\left|c_j-c_{j-1}\right|  \right)^2 \leq \frac{8 t^2}{m} \sum\limits_{n=1}^{\infty}
\left(\frac{1}{n} \sum\limits_{j=1}^{n}\left|c_j-c_{j-1}\right| \right)^2.
\end{align*}

The key step in the proof is the application here of the discrete Hardy inequality~(\ref{1}) for $p=2$,
\begin{equation}\label{Hardy p=2}
    \sum\limits_{n=1}^{\infty} \left(\frac{1}{n} \sum\limits_{j=1}^{n} a_j\right)^2 \leq 4 \sum\limits_{n=1}^{\infty} a_n^2.
\end{equation}
So we obtain
\begin{align*}
\Xi_2^2(t) & \leq \frac{8 t^2}{m} \sum\limits_{n=1}^{\infty}
\left(\frac{1}{n} \sum\limits_{j=1}^{n}\left|c_j-c_{j-1}\right| \right)^2 \leq \frac{32 t^2}{m} \sum\limits_{n=1}^{\infty} |c_{n}-c_{n-1}|^2 \\
&\leq \frac{32 M t^2}{m}\left\| \sum\limits_{n=1}^{\infty} (c_{n}-c_{n-1}) e_n  \right\|^2=  \frac{32 M t^2}{m} \|x\|_1^2.
\end{align*}
Combining this estimate with~(\ref{Ksi1}) we arrive at the following,
$$\left\|(\mathfrak{f})\sum\limits_{n=1}^{\infty} e^{i t \ln n} c_{n} e_n\right\|_{1}^2 \leq\left(\Xi_1 + \Xi_2(t) \right)^2 \leq 2 \Xi_1^2 +2\Xi_2^2(t)\leq \frac{2M }{m} \left(1+ 32 t^2\right) \|x\|_1^2.$$

The last estimate shows that we can define a one-parameter family of operators
$e^{A t}\in[H_1\left(\{e_n\}\right)]$, $t\in\mathbb{R}$, by the formula
$$e^{A t} x=e^{A t} (\mathfrak{f})\sum\limits_{n=1}^{\infty} c_n e_n=(\mathfrak{f})\sum\limits_{n=1}^{\infty} e^{i t \ln n} c_{n} e_n,$$
and that for all $t\in\mathbb{R}$ we have the following estimate,
\begin{equation}\label{bounded ubove}
    \left\|e^{A t}\right\| \leq \sqrt{\frac{2M}{m}} \sqrt{1+ 32 t^2}.
\end{equation}

The next step we prove that one-parameter family of operators $\{e^{A t}\}_{t\in \mathbb{R}}$ is strongly continuous at zero.
For this purpose we observe that for each $x=(\mathfrak{f})\sum\limits_{n=1}^{\infty} c_n e_n\in H_1\left(\{e_n\}\right)$ we have
\begin{align*}
&\|e^{A t}x-x\|_1 =\left\|\sum\limits_{n=2}^{\infty} ((e^{i t \ln n}-1) c_{n}-(e^{i t \ln {(n-1)}}-1) c_{n-1}) e_n\right\|\\
&=\Biggl\|\sum\limits_{n=2}^{\infty} ((e^{i t \ln n}-1) (c_{n}-c_{n-1})e_n+\sum\limits_{n=2}^{\infty}(e^{i t \ln n}-e^{i t \ln {(n-1)}}) c_{n-1} e_n \Biggr\|\\
&\leq \frac{1}{\sqrt{m}}\left(\sum\limits_{n=1}^{\infty} |e^{i t \ln n}-1|^2 |c_n-c_{n-1}|^2\right)^{\frac{1}{2}} +\left\|\sum\limits_{n=2}^{\infty} (e^{i t \ln n}-e^{i t \ln (n-1)}) c_{n-1} e_n\right\|\\
& = \Upsilon(t)+ \Xi_2(t).
\end{align*}
Consider the operator $G: H \supset D(G)\mapsto H,$ defined as
$$Gy = G \sum\limits_{n=1}^{\infty} a_n  e_n= \sum\limits_{n=1}^{\infty} i\ln n\cdot a_n e_n,\quad y=\sum\limits_{n=1}^{\infty} a_n  e_n\in H,
$$
with domain
$D(G)=\left\{y= \sum\limits_{n=1}^{\infty} a_n e_n \in H:\:\{ a_n \ln n\}_{n=1}^{\infty}\in \ell_2\right\}.$
It is not hard to show that $G$ generates $C_0$-group $\{T(t)\}_{t\in \mathbb{R}}$ on $H$, which acts by the formula
$$T(t)y=\sum\limits_{n=1}^{\infty} e^{i t \ln n} a_n e_n,\quad y=\sum\limits_{n=1}^{\infty} a_n  e_n\in H,
$$
see also~\cite{Curtain}. Consequently,
\begin{equation}\label{infinites}
 \sum\limits_{n=2}^{\infty} \left|\frac{e^{i t \ln n}-1}{t} - i \ln n \right|^2 |a_n|^2 \rightarrow 0
\end{equation}
and $\sum\limits_{n=1}^{\infty} |e^{i t \ln n}-1|^2 |a_n|^2 \rightarrow 0$, when $t\rightarrow 0,$
hence
$\Upsilon(t) \rightarrow 0,$ when $t\rightarrow 0.$
Since $\Xi_2(t)\rightarrow 0$, when $t\rightarrow 0,$ then the property of strong continuity at zero for the family $\{e^{A t}\}_{t\in \mathbb{R}}$ is proved.
Note that $e^{A \cdot0}=I$ and $\{e^{A t}\}_{t\in \mathbb{R}}$ obviously satisfy the group property.
Thus $\{e^{A t}\}_{t\in \mathbb{R}}$ is a $C_0$-group on $H_1\left(\{e_n\}\right)$.

Let us check that $A$ is an infinitesimal generator of $C_0$-group $\{e^{A t}\}_{t\in \mathbb{R}}$ constructed above.
To prove this fact we have to verify that $D(A)=E,$
where $E=\left\{x\in H_1\left(\{e_n\}\right): \exists \lim\limits_{t\rightarrow 0} \frac{e^{A t}x - x}{t} \right\}$, and
that for each $x\in D(A)$ we have
\begin{equation}\label{Ax}
Ax=\lim\limits_{t\rightarrow 0} \frac{e^{A t}x - x}{t}.
\end{equation}

Let $x=(\mathfrak{f})\sum\limits_{n=1}^{\infty} c_n e_n\in D(A)$. Then, by~(\ref{Domain_A}) we have
$$\sum\limits_{n=2}^{\infty} \left|c_n \ln n - c_{n-1} \ln(n-1) \right|^2 < \infty.$$
Denote $\gamma_n=c_n \ln n$, $n\in \mathbb{N}.$ Then $\{\gamma_n\}_{n=1}^{\infty}\in \ell_2(\Delta)$ and $c_n=\frac{\gamma_n}{\ln n},$ $n\geq 2.$
Consequently, for all $n\geq 3$ we will have
\begin{equation}\label{gamma}
 \left(\ln n + \ln(n-1) \right) |c_{n-1}|=\frac{\left|\gamma_{n-1}\right|  \ln n}{\ln(n-1)} + \frac{\left|\gamma_{n-1}\right| \ln(n-1)}{\ln(n-1)}   \leq 3 \left|\gamma_{n-1}\right|.
\end{equation}

The next step we prove that~(\ref{Ax}) is true:
\begin{align*}
&\left\|\frac{e^{A t}x - x}{t} - Ax \right\|_1^2 = \left\|(\mathfrak{f})\sum\limits_{n=1}^{\infty}  \left(\frac{e^{i t \ln n}-1}{t} -i \ln n\right)c_n e_n \right\|_1^2 \\
&\leq \frac{1}{m} \sum\limits_{n=2}^{\infty} \left| \left(\frac{e^{i t \ln n}-1}{t} -i \ln n\right)c_n - \left(\frac{e^{i t \ln (n-1)}-1}{t} -i \ln (n-1)\right)c_{n-1} \right|^2\\
&\leq \frac{2}{m} \sum\limits_{n=2}^{\infty} \left|\frac{e^{i t \ln n}-1}{t} - i \ln n \right|^2 |c_n-c_{n-1}|^2\\
&+\frac{2}{m} \sum\limits_{n=2}^{\infty} \left| \frac{e^{i t \ln n} - e^{i t \ln (n-1)}}{t} - \left(i\ln n - i \ln (n-1) \right)\right|^2 |c_{n-1}|^2=\Xi_3(t) +  \frac{2\Theta(t)}{m}.
\end{align*}
We note that, by using~(\ref{infinites}), we have that $\Xi_3(t)\rightarrow  0$ as $t\rightarrow 0$.
Therefore
it is sufficient to prove that $\Theta(t)\rightarrow 0$ as $t\rightarrow  0.$

To this end we observe that for all $n\geq 2$ we have
\begin{align*}
&\left| \frac{e^{i t \ln n} - e^{i t \ln (n-1)}}{t} - i\ln n + i \ln (n-1) \right|^2=\frac{\left(\cos\left(t \ln n\right)-\cos\left(t \ln (n-1)\right) \right)^2}{t^2}\\
&+\left(\frac{\sin\left(t \ln n\right)-\sin\left(t \ln (n-1)\right) }{t}  -\ln n + \ln(n-1)\right)^2= \widetilde{\theta}_n(t)+\widetilde{\widetilde{\theta}}_n(t).
\end{align*}
Further we estimate functions $\widetilde{\theta}_n(t)$ for all $t\in \mathbb{R}$ and $n\geq 2$ in a following manner,
\begin{align*}
\widetilde{\theta}_n(t) &=\frac{4}{t^2} \left(\sin \frac{t\left(\ln n + \ln (n-1) \right)}{2}\right)^2 \left(\sin \frac{t\left(\ln n - \ln (n-1) \right)}{2}\right)^2 \\
&\leq 2^{-2} \left(\ln n + \ln (n-1) \right)^2 \left(\ln n - \ln (n-1) \right)^2 t^2.
\end{align*}
By the mean value theorem we conclude that
functions $\widetilde{\widetilde{\theta}}_n(t)$ are estimated for all $t\in \mathbb{R}$ and $n\geq 2$ in a following manner,
\begin{align*}
\widetilde{\widetilde{\theta}}_n(t) &= \left(\frac{\left(\sin\left(t \ln n\right)-t\ln n\right) -  \left(\sin\left(t \ln (n-1)\right) - t\ln(n-1) \right)}{t}\right)^2 \\
&\leq \frac{\left(\cos\left(\theta t\ln n + (1-\theta) t \ln(n-1) \right)-1 \right)^2   \left(t \ln n - t \ln (n-1) \right)^2}{t^2}\\
&\leq  \left(\theta t\ln n + (1-\theta) t \ln(n-1) \right)^2 \left(\ln n - \ln (n-1) \right)^2\\
&\leq  \left(\ln n +  \ln (n-1) \right)^2 \left(\ln n - \ln (n-1) \right)^2 t^2,\quad \text{where}\quad \theta \in \left(0,1 \right).
\end{align*}

These estimates, by using~(\ref{gamma}), taking into account that $n^2 \left(\ln n - \ln (n-1) \right)^2 \leq 4$ for all $n\geq 2$,
lead us to the following,
\begin{align*}
\Theta(t) &=\sum\limits_{n=2}^{\infty} \left( \widetilde{\theta}_n(t)+ \widetilde{\widetilde{\theta}}_n(t)\right) |c_{n-1}|^2\\
&\leq \frac{5 t^2 }{4} \left( \left(\ln2\right)^4 |c_1|^2 + \sum\limits_{n=3}^{\infty} \left(\ln n + \ln (n-1) \right)^2 |c_{n-1}|^2 \left(\ln n - \ln (n-1) \right)^2\right) \\
&\leq \frac{45 t^2 }{4} \left( \left(\ln2\right)^4 |c_1|^2 +  \sum\limits_{n=3}^{\infty} n^2 \left(\ln n - \ln (n-1) \right)^2 \frac{|\gamma_{n-1}|^2}{n^2}\right)\\
&\leq 45 t^2 \left( \left(\ln2\right)^4 |c_1|^2 +  \sum\limits_{n=3}^{\infty} \frac{|\gamma_{n-1}|^2}{n^2}\right).
\end{align*}
Using arguments which are similar to those were used while estimating $\Xi_2^2(t)$ we conclude that
$$\Theta(t) \leq 180 t^2 \left( \left(\ln2\right)^4 |c_1|^2 +  \sum\limits_{n=2}^{\infty} |\gamma_n-\gamma_{n-1}|^2\right) \rightarrow 0$$
as $t\rightarrow 0,$
since $\{\gamma_n\}_{n=1}^{\infty}\in \ell_2(\Delta).$ Thereby $\Theta(t)\rightarrow 0 $ as $t\rightarrow 0,$ $D(A)\subset E$ and $Ax=\lim\limits_{t\rightarrow  0} \frac{e^{A t}x - x}{t}$, when $x\in D(A).$

Conversely, let $x=(\mathfrak{f})\sum\limits_{n=1}^{\infty} c_n e_n\in E$. Denote $z=\lim\limits_{t\rightarrow  0} \frac{e^{A t}x - x}{t} \in H_1\left(\{e_n\}\right).$ Then
$z=(\mathfrak{f})\sum\limits_{n=1}^{\infty} z_n e_n,\quad\text{where} \quad\{z_n\}_{n=1}^{\infty}\in \ell_2(\Delta).$
The condition $\left\| \frac{e^{A t}x - x}{t} -z\right\|_1 \rightarrow 0,$ $t\rightarrow  0$ implies
$$|z_1|^2+\sum\limits_{n=2}^{\infty} \left|\frac{e^{i t \ln n}-1}{t}c_n - z_n - \frac{e^{i t \ln (n-1)}-1}{t}c_{n-1}+z_{n-1} \right|^2\rightarrow 0,\quad t\rightarrow  0.$$
Hence, $z_1=0$ and for all $n\geq 2$ we have
$$\left|\frac{e^{i t \ln n}-1}{t}c_n - z_n - \frac{e^{i t \ln (n-1)}-1}{t}c_{n-1}+z_{n-1} \right|\rightarrow 0,\quad t\rightarrow 0.$$
Remark that, for $n\geq 2$,
$$\left|\frac{e^{i t \ln n}-1}{t}c_n - z_n\right| \leq \left| \frac{e^{i t \ln n}-1}{t}c_n - z_n - \frac{e^{i t \ln (n-1)}-1}{t}c_{n-1}+z_{n-1}\right|
+ \left| \frac{e^{i t \ln (n-1)}-1}{t}c_{n-1} - z_{n-1}\right|.$$ Since $ \frac{e^{i t \ln 1}-1}{t}c_{1} - z_{1}=0,$
for $n=2$ we get from the last inequality that $\left|\frac{e^{i t \ln 2}-1}{t}c_2 - z_2\right|\rightarrow 0$ as $t\rightarrow 0$.
Applying this inequality subsequently for $n=3,4,\ldots$, we finally get the relations
$\left|\frac{e^{i t \ln n}-1}{t}c_n - z_n\right|\rightarrow 0,$ $t\rightarrow 0,$ $n\geq 2$.
By passing to the limit as $t\rightarrow 0$ we obtain
$$z_n=i c_n \ln n,\quad n\in \mathbb{N}.$$
Therefore $z=(\mathfrak{f})\sum\limits_{n=1}^{\infty} i \ln n \cdot c_n  e_n$. It means that $x\in D(A)$ and $z=Ax.$

Thus $A$ is the generator of constructed $C_0$-group $\{e^{A t}\}_{t\in \mathbb{R}}$ and the theorem is completely proved. \hfill$\Box$
\end{pf}

Concerning Theorem~\ref{not-basis-ln} we note the following. It turns out that, even if we consider the spectrum $\{\lambda_n\}_{n=1}^{\infty}$ of $A$ defined by (\ref{A},\ref{Domain_A}) of the same geometric nature, i.e. satisfying
\begin{equation}\label{geometric cond}
  \lim\limits_{n\rightarrow \infty} i \lambda_n=-\infty\quad\text{and}\quad\lim\limits_{n\rightarrow \infty} |\lambda_{n+1}-\lambda_n|=0,
\end{equation}
then $A$ not necessary generates even a $C_0$-semigroup on $H_1\left(\{e_n\}\right)$. To show this we choose $\lambda_n=i \sqrt{n}$ and prove the following.
\begin{prop}\label{sqrt}
The operator $A$ defined by (\ref{A}) with domain (\ref{Domain_A}), where $\lambda_n=i \sqrt{n}$, $n\in\mathbb{N},$
does not generate the $C_0$-semigroup on $H_1\left(\{e_n\}\right)$.
\end{prop}
\begin{pf}
Suppose that $A$ generate the $C_0$-semigroup $\{S(t)\}_{t\geq 0}$ on $H_1\left(\{e_n\}\right)$.
Since $A e_n=i \sqrt{n} \cdot e_n,\: n\in\mathbb{N},$ for every $t\geq 0$ we will have $S(t) e_n = e^{i t \sqrt{n}} e_n,\:n\in\mathbb{N}.$
Consequently,
$$S(t)\sum\limits_{n=1}^{N} c_n e_n = \sum\limits_{n=1}^{N} e^{i t \sqrt{n}} c_n e_n,\quad N\in\mathbb{N}.$$

Next we show that $S(1)$ is an unbounded operator. To this end we note that, by~(\ref{Riesz}) and using the triangle inequality, the following holds:
\begin{align*}
\left\| S(1)\sum\limits_{n=1}^{N} c_n e_n\right\|_1 &= \left\|\sum\limits_{n=1}^{N+1} \left( e^{i \sqrt{n}} c_n - e^{i \sqrt{n-1}} c_{n-1}\right) e_n \right\|\\
&\geq \left\| \sum\limits_{n=2}^{N+1} (e^{i \sqrt{n}}-e^{i \sqrt{n-1}}) c_{n-1} e_n\right\| - \left\| \sum\limits_{n=1}^{N+1} e^{i \sqrt{n}} (c_n-c_{n-1}) e_n\right\|\\
&\geq \Xi_N -\sqrt{\frac{M}{m}} \left\|\sum\limits_{n=1}^{N} c_n e_n\right\|_1.
\end{align*}

Further, since $\sin s\geq s/2$ for $s\in[0,1]$, and $n \left(\sqrt{n} -\sqrt{n-1} \right)^2 \geq \frac{1}{4}$ for $n\geq 2,$ then
\begin{align*}
\Xi_N^2 &\geq \frac{1}{M} \sum\limits_{n=2}^{N+1} \left| 1-e^{i \left(\sqrt{n-1}-\sqrt{n}\right)}\right|^2 \left|c_{n-1}\right|^2
\geq \frac{1}{M} \sum\limits_{n=2}^{N+1} \sin^2 \left(\sqrt{n-1}-\sqrt{n}\right) \left|c_{n-1}\right|^2\\
&\geq \frac{1}{4M}\sum\limits_{n=2}^{N+1}  n \left(\sqrt{n} -\sqrt{n-1} \right)^2 \frac{\left|c_{n-1}\right|^2}{n}\geq \frac{1}{16 M}\sum\limits_{n=2}^{N+1}
\frac{\left|c_{n-1}\right|^2}{n}.
\end{align*}
Now choose the sequence $x^N=\sum\limits_{n=1}^{N} e_n$, $N\in\mathbb{N}$. Then we see that, on the one hand, $\sup\limits_{N} \left\|x^N \right\|_1 \leq \sqrt{\frac{2}{m}},$
and, on the other hand,
$$
\left\| S(1)x^N\right\|_1=\left\| S(1)\sum\limits_{n=1}^{N} e_n\right\|_1 \geq \frac{1}{4}\sqrt{\frac{1}{ M}\sum\limits_{n=2}^{N+1}
\frac{1}{n}} -\frac{\sqrt{2M}}{m} \rightarrow \infty,\quad N\rightarrow \infty.
$$
Thus we arrive at the unboundedness of $S(1)$, and hence to a contradiction with the definition of the $C_0$-semigroup.
Consequently  $A$ does not generate a $C_0$-semigroup on $H_1\left(\{e_n\}\right)$. \hfill$\Box$
\end{pf}

The underlying reason of the phenomenon described above is as follows. The rotations $e^{i t \sqrt{x}},$ are slowing down, when $x\rightarrow +\infty$, with too low speed to guarantee the convergence of the series
$$\sum\limits_{n=2}^{\infty} (e^{i t \sqrt{n}}-e^{i t \sqrt{n-1}}) c_{n-1} e_n$$
in $H$. Using the same arguments as in the proof of Proposition~\ref{sqrt} we can say the following.
\begin{prop}\label{sqrt h}
Let $\{\lambda_n\}_{n=1}^{\infty}\subset i\mathbb{R}$
satisfy~(\ref{geometric cond}) and suppose that there exists $\alpha\in \left(0,\frac{1}{2}\right]$ such that $\liminf\limits_{n\rightarrow \infty} n^{\alpha} |\lambda_n - \lambda_{n-1}|>0.$
Then the operator $A$ defined by (\ref{A}) with domain (\ref{Domain_A})
does not generate a $C_0$-semigroup on $H_1\left(\{e_n\}\right)$.
\end{prop}

\subsection{Generalization of the Theorem~\ref{not-basis-ln} and infinitesimal operators on $H_k\left(\{e_n\}\right)$}
We generalize the Theorem~\ref{not-basis-ln} in two directions. First, we give some other constructions of infinitesimal operators with complete minimal non-basis family of eigenvectors on $H_k\left(\{e_n\}\right)$ spaces.
Second, we consider more general behaviour of the spectrum of these operators. Throughout the paper it is considered that, if $\left\{u_n\right\}_{n=1}^{\infty}$ is a real sequence,
then $u_{-n}=0$ for $n\geq 0.$

Recall that
$$\Delta^k u_n= \sum\limits_{j=0}^{k} (-1)^j C_{k}^j u_{n-j},\quad k\in\mathbb{Z}_{+}.$$
For the generalization of the Theorem~\ref{not-basis-ln} we introduce the following classes of sequences.

\begin{definition}\label{S_k} Let $k\in\mathbb{N}$ and let $\left\{f(n)\right\}_{n=1}^{\infty}\subset \mathbb{R}$ be a sequence.
Then we define
$$\mathcal{S}_k=\Bigl\{\left\{f(n)\right\}_{n=1}^{\infty}:\: \lim\limits_{n\rightarrow \infty} f(n)=+\infty; \left\{n^j \Delta^j f(n) \right\}_{n=1}^{\infty}\in \ell_{\infty}\:\: \text{for} \:\: 1\leq j\leq k\Bigr\}.$$
\end{definition}

It is clear that $\mathcal{S}_k\subseteq \mathcal{S}_m$ provided $m\leq k$.
For example, for each $k\in\mathbb{N}$ we have $\left\{\ln n\right\}_{n=1}^{\infty}\in \mathcal{S}_k$, $\left\{\ln\ln(n+1)\right\}_{n=1}^{\infty}\in \mathcal{S}_k,$ $\left\{\ln\ln\sqrt{n+1}\right\}_{n=1}^{\infty}\in \mathcal{S}_k$, and $\left\{\sqrt{n}\right\}_{n=1}^{\infty}\notin \mathcal{S}_k$.

Now we formulate our generalization.
\begin{thm}\label{not-basis_fk}
Let $k\in\mathbb{N}$. Then the operator $A_k:H_k\left(\{e_n\}\right) \supset D(A_k) \mapsto H_k\left(\{e_n\}\right),$ defined by
$$A_k x=A_k (\mathfrak{f})\sum\limits_{n=1}^{\infty} c_{n} e_n= (\mathfrak{f})\sum\limits_{n=1}^{\infty} i f(n) \cdot c_{n} e_n,$$
where $\left\{f(n)\right\}_{n=1}^{\infty}\in\mathcal{S}_k$, with domain
\begin{equation}\label{Domain_k}
    D(A_k)=\left\{x= (\mathfrak{f})\sum\limits_{n=1}^{\infty} c_{n} e_n \in H_k\left(\{e_n\}\right):\:  \{f(n) \cdot c_{n}\}_{n=1}^{\infty}\in \ell_2(\Delta^k)\right\},
\end{equation}
generates the $C_0$-group on $H_k\left(\{e_n\}\right)$, which acts for every $t\in\mathbb{R}$ by the formula
$$e^{A_k t} x=e^{A_k t}(\mathfrak{f})\sum\limits_{n=1}^{\infty} c_n e_n =(\mathfrak{f})\sum\limits_{n=1}^{\infty} e^{i t f(n)} c_{n} e_n.
$$
\end{thm}
\begin{pf} For the case when $k=1$ the proof is similar to the proof of Theorem~\ref{not-basis-ln}.

Let us prove the theorem for any fixed $k\geq2.$
Consider
$x=(\mathfrak{f})\sum\limits_{n=1}^{\infty} c_n e_n\in H_k\left(\{e_n\}\right)$. Fix $t\in\mathbb{R}$ and prove that
$(\mathfrak{f})\sum\limits_{n=1}^{\infty} e^{i t f(n)} c_{n} e_n \in H_k\left(\{e_n\}\right).$
By~(\ref{norm}) we have
\begin{align*}
\left\|(\mathfrak{f})\sum\limits_{n=1}^{\infty} e^{i t f(n)} c_{n} e_n\right\|_{k} &= \Biggl\|\sum\limits_{n=1}^{\infty}\left(\Delta^k \left(e^{i t f(n)}c_{n}\right)\right) e_n\Biggr\|,
\end{align*}
where $c_{1-j}=0$ for any $j\in\mathbb{N}$.
Further we will use the following Leibnitz theorem for finite differences:
\begin{equation}\label{Leibnitz}
\Delta^k (u_n v_n) = \sum\limits_{j=0}^{k} C_{k}^j \Delta^{k-j} u_{n-j} \Delta^j v_{n},\quad k\in \mathbb{N},
\end{equation}
see.~\cite{MTLM}, pp. 34-35.
By Leibnitz theorem~(\ref{Leibnitz}), taking into account~(\ref{Riesz}), since $c_{n-j}=0$ for $j\geq n$, we have
\begin{align*}
\left\|(\mathfrak{f})\sum\limits_{n=1}^{\infty} e^{i t f(n)} c_{n} e_n\right\|_{k} &\leq \sum\limits_{j=0}^{k} C_{k}^j  \left\| \sum\limits_{n=1}^{\infty} \left(\Delta^j e^{i t f(n)}\right) \left(\Delta^{k-j} c_{n-j} \right) e_n \right\|\\
&\leq \frac{1}{\sqrt{m}} \sum\limits_{j=0}^{k} C_{k}^j \left(\sum\limits_{n=1}^{\infty} \left|\Delta^j e^{i t f(n)}\right|^2 \left|\Delta^{k-j} c_{n-j} \right|^2 \right)^{\frac{1}{2}}\\
&= \frac{1}{\sqrt{m}} \sum\limits_{j=0}^{k} C_{k}^j \left(\sum\limits_{n=j+1}^{\infty} \left|\Delta^j e^{i t f(n)}\right|^2 \left|\Delta^{k-j} c_{n-j} \right|^2 \right)^{\frac{1}{2}}\\
&= \frac{1}{\sqrt{m}} \sum\limits_{j=0}^{k} C_{k}^j \sqrt{\Theta_j(t)}.
\end{align*}
Function $\Theta_0(t)$ is estimated by~(\ref{Riesz}) and~(\ref{norm}) immediately:
$$\Theta_0(t) = \sum\limits_{n=1}^{\infty} \left|e^{i t f(n)} \Delta^{k} c_{n} \right|^2=  \sum\limits_{n=1}^{\infty} \left| \Delta^{k} c_{n} \right|^2\leq M \left\| \sum\limits_{n=1}^{\infty} \Delta^{k} c_{n} e_n \right\|^2 =M \|x\|_k^2.
$$

Further we need to estimate functions
$\Theta_j(t)$, $1\leq j\leq k$. For this purpose we note that, since $\left\{f(n)\right\}_{n=1}^{\infty}\in\mathcal{S}_k$, then there exists $C> 0$ such that for any fixed $j:$
$1\leq j \leq k$ and
for all $n\in\mathbb{N}$ we have
\begin{equation}\label{nj}
n^j \left|\Delta^j f(n) \right|\leq C.
\end{equation}

Next, since
$\left\{f(n)\right\}_{n=1}^{\infty}\in\mathcal{S}_k$, then for every $m\in\mathbb{N}$ we will have
$\left\{f(n-m)\right\}_{n=1}^{\infty}\in\mathcal{S}_k.$ Indeed, by~(\ref{nj}), taking into account that $f(n-m)=0$ for $1\leq n\leq m$, for all $n\in\mathbb{N}$ we have $(n-m)^j\left|\Delta^j f(n-m) \right|\leq C,$ hence for $n>m$ we get $\left|\Delta^j f(n-m) \right|\leq \frac{C}{(n-m)^j}$.
Note that for $n>m$
$$\frac{1}{(n-m)^j}\leq \frac{(1+m)^j}{n^j}.$$
It follows from the following fact:
$$\left(\frac{n}{n-m} \right)^j=\left(1+\frac{m}{n-m} \right)^j\leq \left(1+\frac{m}{m+1-m} \right)^j=(1+m)^j.$$

Thus for $n>m$ we obtain
$\left|\Delta^j f(n-m) \right|\leq \frac{C(1+m)^j}{n^j}$, and for $1\leq n\leq m$ we will get $\left|\Delta^j f(n-m) \right|=0$. Consequently,
$n^j\left|\Delta^j f(n-m) \right|\leq C(1+m)^j$ for all
$n$, hence $\left\{f(n-m)\right\}_{n=1}^{\infty}\in\mathcal{S}_k.$

Let  $m:$ $1\leq m\leq k$.
Consider the following sets:
$\Sigma_1=\{0,1,2,\dots,k-1\}$, $\Sigma_2=\{0,1,2,\dots,k-2 \},\dots,$ $\Sigma_{k-1}=\{0,1\},$ $\Sigma_k=\{0\}$. It is clear that
$\Sigma_1\supset \Sigma_2 \supset \Sigma_3 \supset \dots \supset \Sigma_{k}.$

Now we claim that for every $m:$ $1\leq m\leq k$ and every sequence $\left\{\widetilde{f}(n)\right\}_{n=1}^{\infty}\in\mathcal{S}_k$,
for any $s\in \Sigma_m$, $t\in\mathbb{R}$ and arbitrary $n>m$ the following estimate holds true:
\begin{equation}\label{Exponents S}
\left| \Delta^m e^{(-1)^{s}i t \Delta^{s} \widetilde{f}(n)}\right| \leq \frac{\mathcal{P}_m\left[\widetilde{f}(n)\right](|t|)}{n^{m+s}},
\end{equation}
where $\mathcal{P}_m\left[\widetilde{f}(n)\right]$ is a polynomial of degree $m$, with positive coefficients depending on $\left\{\widetilde{f}(n)\right\}_{n=1}^{\infty}$, and without a free term.
We will prove this fact by induction.

Further for convenience we will denote $s$ by $s(m),$ if $s\in \Sigma_m$.

The induction base. Let $m=1.$ In this case we have a set $\Sigma_1=\{0,1,2,3,\dots,k-1\}$ and $s=s(1)$ takes values from $0$ to $k-1$. Then for every $\left\{\widetilde{f}(n)\right\}_{n=1}^{\infty}\in\mathcal{S}_k$
and any $s(1)\in \Sigma_1$, for $n>1$, taking into account~(\ref{nj}), we will have
\begin{align*}
&\left| \Delta e^{(-1)^{s(1)}i t \Delta^{s(1)} \widetilde{f}(n)}\right|= \left| e^{(-1)^{s(1)}i t \Delta^{s(1)} \widetilde{f}(n)} -  e^{(-1)^{s(1)}i t \Delta^{s(1)} \widetilde{f}(n-1)} \right|\\
& = \left|e^{(-1)^{s(1)}i t \Delta^{s(1)} \widetilde{f}(n)} \left(1-e^{(-1)^{s(1)+1}i t \Delta^{s(1)+1} \widetilde{f}(n)} \right) \right| = \left|1-e^{(-1)^{s(1)+1}i t \Delta^{s(1)+1} \widetilde{f}(n)} \right|\\
&\leq \sqrt{2} |t| \left| \Delta^{s(1)+1} \widetilde{f}(n)\right| = \frac{\sqrt{2} |t| n^{s(1)+1} \left| \Delta^{s(1)+1} \widetilde{f}(n)\right|}{n^{s(1)+1}} \leq \frac{\sqrt{2} C\left(\widetilde{f}\right) |t|}{n^{s(1)+1}},
\end{align*}
where $C\left(\widetilde{f}\right)$ is a constant depending on the sequence $\left\{\widetilde{f}(n)\right\}_{n=1}^{\infty}.$

Now assume that for arbitrary
$l:$ $1\leq l \leq m-1$ and every  $\left\{\widetilde{f}(n)\right\}_{n=1}^{\infty}\in\mathcal{S}_k$, for all $s(l)\in \Sigma_l$, $t\in\mathbb{R}$ and
$n>l$ the following is valid:
\begin{equation}\label{Exponents S assumption}
\left| \Delta^l e^{(-1)^{s(l)}i t \Delta^{s(l)} \widetilde{f}(n)}\right| \leq \frac{\mathcal{P}_l\left[\widetilde{f}(n)\right](|t|)}{n^{l+s(l)}},
\end{equation}
where $\mathcal{P}_l\left[\widetilde{f}(n)\right]$  is a polynomial of degree $l$, with positive coefficients depending on $\left\{\widetilde{f}(n)\right\}_{n=1}^{\infty}$, and without a free term.

Next we will prove (\ref{Exponents S assumption}) for $l=m.$ For any $\left\{\widetilde{f}(n)\right\}_{n=1}^{\infty}\in\mathcal{S}_k$, arbitrary $s(m)\in \Sigma_m$ and for $n>m$, by Leibnitz theorem~(\ref{Leibnitz}) we have that
\begin{align*}
&\left| \Delta^m e^{(-1)^{s(m)}i t \Delta^{s(m)} \widetilde{f}(n)}\right|= \left|\Delta^{m-1} \left(e^{(-1)^{s(m)}i t \Delta^{s(m)} \widetilde{f}(n)} -e^{(-1)^{s(m)}i t \Delta^{s(m)} \widetilde{f}(n-1)} \right) \right|\\
& = \left|\Delta^{m-1} \left( e^{(-1)^{s(m)}i t \Delta^{s(m)} \widetilde{f}(n)}\left(1-e^{(-1)^{s(m)+1}i t \Delta^{s(m)+1} \widetilde{f}(n)} \right) \right) \right|\\
& = \left| \sum\limits_{l=0}^{m-1} C_{m-1}^{l} \left( \Delta^l e^{(-1)^{s(m)}i t \Delta^{s(m)} \widetilde{f}(n)} \right) \left(\Delta^{m-1-l} \left(1- e^{(-1)^{s(m)+1}i t \Delta^{s(m)+1} \widetilde{f}(n-l)}\right)  \right) \right|\\
&\qquad \leq\left| e^{(-1)^{s(m)}i t \Delta^{s(m)} \widetilde{f}(n)}\right| \left|\Delta^{m-1} e^{(-1)^{s(m)+1}i t \Delta^{s(m)+1} \widetilde{f}(n)} \right|\\
&\qquad\qquad + \sum\limits_{l=1}^{m-2} C_{m-1}^{l} \left| \Delta^l e^{(-1)^{s(m)}i t \Delta^{s(m)} \widetilde{f}(n)}\right| \left|\Delta^{m-1-l} e^{(-1)^{s(m)+1}i t \Delta^{s(m)+1} \widetilde{f}(n-l)} \right|\\
&\qquad\qquad\qquad+ \left| \Delta^{m-1} e^{(-1)^{s(m)}i t \Delta^{s(m)} \widetilde{f}(n)}\right| \left|1- e^{(-1)^{s(m)+1}i t \Delta^{s(m)+1} \widetilde{f}(n-(m-1))} \right|.
\end{align*}

Since $s(m)\in \Sigma_{m},$ and $\Sigma_m \subset \Sigma_{m-1} \subset \Sigma_{m-2} \subset \dots \subset \Sigma_2 \subset \Sigma_1$, then $s(m)\in \Sigma_l$ for any $l:$ $1\leq l\leq m-1.$ Analogously, since $s(m)+1 \in \Sigma_{m-1},$ then $s(m)+1\in \Sigma_l$ for any $l:$ $1\leq l\leq m-1.$  Earlier we proved that, if  $\left\{\widetilde{f}(n)\right\}_{n=1}^{\infty}\in\mathcal{S}_k$, then for every $m\in\mathbb{N}$ we have that
$\left\{\widetilde{f}(n-m)\right\}_{n=1}^{\infty}\in\mathcal{S}_k.$
Thus we can use the induction assumption~(\ref{Exponents S assumption}) and continue the estimate,
\begin{align*}
\left| \Delta^m e^{(-1)^{s(m)}i t \Delta^{s(m)} \widetilde{f}(n)}\right| &\leq  \frac{\mathcal{P}_{m-1}\left[\widetilde{f}(n)\right](|t|)}{n^{m-1+s(m)+1}}+\sum\limits_{l=1}^{m-2} C_{m-1}^{l} \frac{\mathcal{P}_{l}\left[\widetilde{f}(n)\right](|t|)}{n^{l+s(m)}}  \\
&\times\frac{\mathcal{P}_{m-1-l}\left[\widetilde{f}(n-l)\right](|t|)}{n^{m-1-l+s(m)+1}}+\frac{\mathcal{P}_{m-1}\left[\widetilde{f}(n)\right](|t|)}{n^{m-1+s(m)}} \frac{\sqrt{2} \widetilde{C} |t|}{n^{s(m)+1}}.
\end{align*}
Therefore we have that
$$\left| \Delta^m e^{(-1)^{s(m)}i t \Delta^{s(m)} \widetilde{f}(n)}\right|\leq \frac{\widetilde{P_m}\left[\widetilde{f}(n)\right](|t|)}{n^{m+s(m)}},$$
where $\widetilde{P_m}\left[\widetilde{f}(n)\right]$ is a polynomial of degree  $m$, with positive coefficients depending on $\left\{\widetilde{f}(n)\right\}_{n=1}^{\infty}$, and without a free term, hence~(\ref{Exponents S}) is proved.

It follows from~(\ref{Exponents S}), in particular, that for our sequence $\left\{f(n)\right\}_{n=1}^{\infty}\in\mathcal{S}_k$,
for any $j:$ $1\leq j \leq k$, $t\in\mathbb{R}$ and for every $n>j$ the following estimate is true:
\begin{equation}\label{Exponents}
\left| \Delta^j e^{i t f(n)}\right| \leq \frac{\mathcal{P}_j (|t|)}{n^j},
\end{equation}
where $\mathcal{P}_j$ is a polynomial of degree $j$, with positive coefficients and without a free term.

Further, on the basis of~(\ref{Exponents}) functions $\Theta_j(t)$,  $1\leq j \leq k$, are estimated in the following way:
\begin{equation}\label{theta}
\Theta_j(t)\leq \left(\mathcal{P}_j (|t|)\right)^2 \sum\limits_{n=j+1}^{\infty} \frac{\left|\Delta^{k-j} c_{n-j} \right|^2}{n^{2j}}= \left(\mathcal{P}_j (|t|)\right)^2 \Omega_j.
\end{equation}
We observe that for all $d,n\in\mathbb{N}$ the formula
\begin{equation}\label{formula}
\Delta^d c_n= \sum\limits_{m=1}^{n} \Delta^{d+1} c_m
\end{equation}
is valid.
Hence we can estimate $\Omega_j$, $1\leq j \leq k$, in the following manner,
\begin{align*}
\Omega_j &= \sum\limits_{n=j+1}^{\infty} \frac{\left|\Delta^{k-j} c_{n-j} \right|^2}{n^{2j}}\leq  \sum\limits_{n=1}^{\infty} \frac{\left|\Delta^{k-j} c_{n-j} \right|^2}{n^{2j}}= \sum\limits_{n=1}^{\infty} \left(\frac{1}{n^j} \left|\sum\limits_{m=1}^{n-j} \Delta^{k-j+1} c_m \right| \right)^2 \\
&\leq \sum\limits_{n=1}^{\infty} \left(\frac{1}{n^j} \sum\limits_{m=1}^{n-j} \left|\Delta^{k-j+1} c_m \right| \right)^2 \leq \sum\limits_{n=1}^{\infty} \left(\frac{1}{n^j} \sum\limits_{m=1}^{n} \left|\Delta^{k-j+1} c_m \right| \right)^2.
\end{align*}
Applying the discrete Hardy inequality~(\ref{Hardy p=2}) we will have
$$\Omega_j \leq 4 \sum\limits_{n=1}^{\infty} \left(\frac{1}{n^{j-1}} \left|\Delta^{k-j+1} c_n \right| \right)^2.$$

To continue the estimate we will  consistently apply formula~(\ref{formula}) and Hardy inequaltity~(\ref{Hardy p=2}) $j-1$ more times:
\begin{align*}
\Omega_j &\leq 4 \sum\limits_{n=1}^{\infty} \left(\frac{1}{n^{j-1}} \left|\Delta^{k-j+1} c_n \right| \right)^2 = 4 \sum\limits_{n=1}^{\infty} \left(\frac{1}{n^{j-1}} \left|\sum\limits_{m=1}^{n} \Delta^{k-j+2} c_m \right| \right)^2\\
&\qquad\leq 4 \sum\limits_{n=1}^{\infty} \left(\frac{1}{n^{j-1}} \sum\limits_{m=1}^{n} \left|\Delta^{k-j+2} c_m \right| \right)^2 \leq 4^2\sum\limits_{n=1}^{\infty} \left(\frac{1}{n^{j-2}} \left|\Delta^{k-j+2} c_n \right| \right)^2\\
&\qquad\qquad \leq \ldots \leq 4^{j-1} \sum\limits_{n=1}^{\infty} \left(\frac{1}{n} \left|\Delta^{k-1} c_n \right| \right)^2 \leq  4^{j-1} \sum\limits_{n=1}^{\infty} \left(\frac{1}{n} \sum\limits_{m=1}^{n} \left|\Delta^{k} c_m \right| \right)^2 \\
&\qquad\qquad\qquad  \leq 4^j\sum\limits_{n=1}^{\infty}  \left|\Delta^{k} c_n \right|^2 \leq 4^j M \left\| \sum\limits_{n=1}^{\infty} \Delta^{k} c_n e_n \right\|^2 =4^j M \|x\|_k^{2}.
\end{align*}

Combining this estimate with~(\ref{theta}) we arrive at the following,
\begin{align*}
\left\|(\mathfrak{f})\sum\limits_{n=1}^{\infty} e^{i t f(n)} c_{n} e_n\right\|_{k} &\leq \frac{1}{\sqrt{m}} \sum\limits_{j=0}^{k} C_{k}^j \sqrt{\Theta_j(t)}
\end{align*}
\begin{equation}\label{series bound}
\leq \frac{1}{\sqrt{m}}\Biggl(
\sqrt{M} \|x\|_k+ \sum\limits_{j=1}^{k} C_{k}^j\sqrt{\left(\mathcal{P}_j (|t|)\right)^2 \Omega_j} \Biggr) \leq \sqrt{\frac{M}{m}}
\sum\limits_{j=0}^{k} 2^j C_{k}^j  \mathcal{P}_j (|t|) \|x\|_k,
\end{equation}
where $\mathcal{P}_0 (|t|)\equiv 1.$
This estimate shows that we can define a one-parametric family of operators $e^{A_k t}\in[H_k\left(\{e_n\}\right)]$, $t\in\mathbb{R}$, by the formula
$$e^{A_k t} x=e^{A_k t} (\mathfrak{f})\sum\limits_{n=1}^{\infty} c_n e_n=(\mathfrak{f})\sum\limits_{n=1}^{\infty} e^{i t f(n)} c_{n} e_n,$$
and that for all $t\in\mathbb{R}$ we have
\begin{equation}\label{bounded ubove k}
    \left\|e^{A_k t}\right\| \leq \mathfrak{p}_k (|t|),
\end{equation}
where $\mathfrak{p}_k (|t|)=\sqrt{\frac{M}{m}}\sum\limits_{j=0}^{k} 2^j C_{k}^j  \mathcal{P}_j (|t|)$ is a polynomial  with positive coefficients  and $\deg \mathfrak{p}_k =k.$

Next we prove that the family $\{e^{A_k t}\}_{t\in \mathbb{R}}$ is strongly continuous at zero.
To this end we note that for any $j:$ $1\leq j \leq k$ and for every $n>j$ we have $\Delta^j \left(e^{i t f(n)}-1 \right)=\Delta^j \left(e^{i t f(n)}\right)$.
Then  by Leibnitz theorem~(\ref{Leibnitz}), taking into account estimate~(\ref{series bound}), for each $x=(\mathfrak{f})\sum\limits_{n=1}^{\infty} c_n e_n\in H_k\left(\{e_n\}\right)$ we will have the following,
\begin{align*}
\|e^{A_k t}x-x\|_k &=\left\|\sum\limits_{n=1}^{\infty} \left(\Delta^k\left((e^{i t f(n)}-1) c_{n}\right)\right) e_n\right\|\\
&\leq \sum\limits_{j=0}^{k} C_{k}^j  \left\| \sum\limits_{n=1}^{\infty} \left(\Delta^j \left(e^{i t f(n)}-1 \right)\right) \left(\Delta^{k-j} c_{n-j} \right) e_n \right\|\\
&\leq\frac{1}{\sqrt{m}}  \Biggl(\left(\sum\limits_{n=1}^{\infty} \left|e^{i t f(n)}-1\right|^2 \left|\Delta^{k} c_{n} \right|^2 \right)^{\frac{1}{2}}\\
&+\sum\limits_{j=1}^{k} C_{k}^j \left(\sum\limits_{n=j+1}^{\infty} \left|\Delta^j e^{i t f(n)}\right|^2 \left|\Delta^{k-j} c_{n-j} \right|^2 \right)^{\frac{1}{2}}\Biggr)=\Upsilon_k(t)\\
& +\frac{1}{\sqrt{m}} \sum\limits_{j=1}^{k} C_{k}^j \sqrt{\Theta_j(t)}\leq \Upsilon_k(t)+ \sqrt{\frac{M}{m}}
\sum\limits_{j=1}^{k} 2^j C_{k}^j  \mathcal{P}_j (|t|) \|x\|_k.
\end{align*}

The function $\widetilde{\mathfrak{p}}_k (|t|)=\sum\limits_{j=1}^{k} 2^j C_{k}^j  \mathcal{P}_j (|t|)$ is a polynomial of argument $|t|$ with degree $k$, with positive coefficients and without a free term. Thus $\sqrt{\frac{M}{m}}
\widetilde{\mathfrak{p}}_k (|t|)\|x\|_k\rightarrow 0$ for $t\rightarrow 0.$

Next we will show that $\Upsilon_k(t) \rightarrow 0$ as $t\rightarrow 0.$
To this end let us consider the operator $F: H \supset D(F)\mapsto H,$ defined by
$$Fy =F \sum\limits_{n=1}^{\infty} a_n  e_n= \sum\limits_{n=1}^{\infty} if(n)\cdot a_n e_n,\quad y=\sum\limits_{n=1}^{\infty} a_n  e_n\in H,
$$
with domain
$D(F)=\left\{y= \sum\limits_{n=1}^{\infty} a_n e_n \in H:\:\{ a_n f(n)\}_{n=1}^{\infty}\in \ell_2\right\}.$

It is not hard to prove that $F$ generates the $C_0$-group $\{T(t)\}_{t\in \mathbb{R}}$ on $H$, which acts by the formula
$$T(t)y=\sum\limits_{n=1}^{\infty} e^{i t f(n)} a_n e_n,\quad y=\sum\limits_{n=1}^{\infty} a_n  e_n\in H,
$$
see also~\cite{Curtain}. Consequently,
$\sum\limits_{n=1}^{\infty} |e^{i t f(n)}-1|^2 |a_n|^2 \rightarrow 0$ as $t\rightarrow 0.$
Since $\left\{\Delta^{k} c_{n}\right\}_{n=1}^{\infty}\in \ell_2$, then
$\Upsilon_k(t) \rightarrow 0,$ as $t\rightarrow 0.$ Hence $\|e^{A_k t}x-x\|_k\rightarrow 0$  as $t\rightarrow 0,$
and family $\{e^{A_k t}\}_{t\in \mathbb{R}}$ is strongly continuous at zero.
It is obvious that $e^{A_k \cdot0}=I$ and the group property holds for $\{e^{A_k t}\}_{t\in \mathbb{R}}$.
Thus $\{e^{A_k t}\}_{t\in \mathbb{R}}$ is the $C_0$-group on the space $H_k\left(\{e_n\}\right)$.

Further we will show that $A_k$ is an infinitesimal generator of constructed $C_0$-group $\{e^{A_k t}\}_{t\in \mathbb{R}}$.
For this purpose we need to check that $D(A_k)=E_k,$
where $E_k=\left\{x\in H_k\left(\{e_n\}\right): \exists \lim\limits_{t\rightarrow 0} \frac{e^{A_k t}x - x}{t} \right\}$ and
that for each $x\in D(A_k)$ we have
\begin{equation}\label{A_k x}
A_k x=\lim\limits_{t\rightarrow 0} \frac{e^{A_k t}x - x}{t}.
\end{equation}

Let $x=(\mathfrak{f})\sum\limits_{n=1}^{\infty} c_n e_n\in E_k$.
Denote $z=\lim\limits_{t\rightarrow  0} \frac{e^{A_k t}x - x}{t} \in H_k\left(\{e_n\}\right).$ Then
$z=(\mathfrak{f})\sum\limits_{n=1}^{\infty} z_n e_n,$ where $\{z_n\}_{n=1}^{\infty}\in \ell_2(\Delta^k).$
The condition $\left\| \frac{e^{A_k t}x - x}{t} -z\right\|_k \rightarrow 0,$ $t\rightarrow  0$, implies
$$\sum\limits_{n=1}^{\infty} \left|\Delta^k \left(\frac{e^{i t f(n)}-1}{t}c_n - z_n \right) \right|^2\rightarrow 0,\quad t\rightarrow  0,$$
hence for every $n\in\mathbb{N}$  we have that
\begin{equation}\label{lim}
 \left|\Delta^k \left(\frac{e^{i t f(n)}-1}{t}c_n - z_n \right) \right|^2\rightarrow 0,\quad t\rightarrow  0.
\end{equation}
We will use induction on $n$ to prove that $z_n=i c_n f(n),\: n\in \mathbb{N}.$

For $n=1$ from~(\ref{lim}) we get $\frac{e^{i t f(1)}-1}{t}c_1 - z_1 \rightarrow 0,$ if $t\rightarrow  0,$ hence passing to the limit yields $z_1=ic_1 f(1)$.

Let $z_n=i c_n f(n)$ for $1\leq n \leq l$.
Note that for any $n\in \mathbb{N}$ we have that
$$
\Delta^k \left(\frac{e^{i t f(n)}-1}{t}c_n - z_n \right)   = \sum_{j=0}^{k} (-1)^j C_{k}^{j} \left(\frac{e^{i t f(n-j)}-1}{t}c_{n-j} - z_{n-j} \right).
$$
Therefore, for $n=l+1$ we will have
\begin{align*}
&\left|\frac{e^{i t f(l+1)}-1}{t}c_{l+1} - z_{l+1}  \right| = \Biggl|  \sum_{j=0}^{k} (-1)^j C_{k}^{j} \left(\frac{e^{i t f(l+1-j)}-1}{t}c_{l+1-j} - z_{l+1-j} \right) \\
&-\sum_{j=1}^{k} (-1)^j C_{k}^{j} \left(\frac{e^{i t f(l+1-j)}-1}{t}c_{l+1-j} - z_{l+1-j} \right)  \Biggr| \\
&\leq \left|\Delta^k \left(\frac{e^{i t f(l+1)}-1}{t}c_{l+1} - z_{l+1}\right)\right|+ \sum_{j=1}^{k} C_{k}^{j} \left|\frac{e^{i t f(l+1-j)}-1}{t}c_{l+1-j} - z_{l+1-j} \right|.
\end{align*}
Taking into account that $z_{l+1-j}=i c_{l+1-j} f(l+1-j)$, by~(\ref{lim}) we obtain that
$z_{l+1}=i c_{l+1} f(l+1).$
Hence $z_n=i c_n f(n),\: n\in \mathbb{N},$ and $z=(\mathfrak{f})\sum\limits_{n=1}^{\infty} i f(n)  c_n  e_n$. It means that $x\in D(A_k)$ and $z=A_k x.$

Conversely, assume that $x=(\mathfrak{f})\sum\limits_{n=1}^{\infty} c_n e_n\in D(A_k)$  and prove that $x\in E_k$ and~(\ref{A_k x}) holds.
To this end we note that for every $N\in\mathbb{N}$,
\begin{align*}
&\left\| \frac{e^{A_k t}x - x}{t} -A_k x\right\|_k =\left\|\sum\limits_{n=1}^{\infty} \left(\Delta^k\left(\left(\frac{e^{i t f(n)}-1}{t} - i f(n)\right)c_{n}\right)\right) e_n\right\|\\
&\quad \leq \frac{1}{\sqrt{m}} \left(\sum\limits_{n=1}^{\infty} \left|\Delta^k\left(\left(\frac{e^{i t f(n)}-1}{t} - i f(n)\right)c_{n}\right)\right|^2\right)^{\frac{1}{2}} \\
&\quad\quad\leq \frac{1}{\sqrt{m}} \Biggl(\sum\limits_{n=1}^{N} \left|\Delta^k\left(\left(\frac{e^{i t f(n)}-1}{t} - i f(n)\right)c_{n}\right)\right|^2\\
&\quad\quad\quad+ 2  \sum\limits_{n=N+1}^{\infty} \left|\Delta^k\left(\frac{e^{i t f(n)}-1}{t}c_{n}\right)\right|^2 + 2 \sum\limits_{n=N+1}^{\infty} \left|\Delta^k\left(i f(n)c_{n}\right)\right|^2\Biggr)^{\frac{1}{2}}.
\end{align*}

Since $x\in D(A_k)$, then $A_k x \in H_k\left(\{e_n\}\right),$ hence there exists a constant $C\geq 0:$ $\left\| A_k x\right\|_k \leq C.$ Consequently
$$\frac{1}{M} \sum\limits_{n=1}^{\infty} \left|\Delta^k\left(i f(n)c_{n}\right)\right|^2 \leq\left\| \sum\limits_{n=1}^{\infty} \left(\Delta^k\left(i f(n)c_{n}\right)\right) e_n\right\|^2=\left\| A_k x\right\|_k^2 \leq C^2,$$
hence
\begin{equation}\label{first sum}
\sum\limits_{n=1}^{\infty} \left|\Delta^k\left(i f(n)c_{n}\right)\right|^2 \leq MC^2.
\end{equation}

Now we observe that for each $t\in\mathbb{R}$ and arbitrary $N\in\mathbb{N}$ the following holds,
\begin{align*}
&\left(\sum\limits_{n=N+1}^{\infty} \left|\Delta^k\left(\frac{e^{i t f(n)}-1}{t}c_{n}\right)\right|^2\right)^{\frac{1}{2}}  \leq \left(\sum\limits_{n=N+1}^{\infty} \left|\frac{e^{i t f(n)}-1}{t}\right|^2 \left|\Delta^{k} c_{n} \right|^2 \right)^{\frac{1}{2}}\\
& +  \sum\limits_{j=1}^{k} 2^j C_{k}^j  \frac{\mathcal{P}_j (|t|)}{|t|}\left( \sum\limits_{n=N+1}^{\infty} \left|\Delta^k c_n \right|^2\right)^{\frac{1}{2}} =
\widehat{\Upsilon}_k^N(t) + \mathcal{Q}_k (|t|) \left( \sum\limits_{n=N+1}^{\infty} \left|\Delta^k c_n \right|^2\right)^{\frac{1}{2}},
\end{align*}
where $c_{n}=0$ for $n\leq N$, and $\mathcal{Q}_k(|t|) =  \sum\limits_{j=1}^{k} 2^j C_{k}^j  \frac{\mathcal{P}_j (|t|)}{|t|} $ is a polynomial of argument $|t|$ with degree $k-1$
(see the proof of strong continuity of $\{e^{A_k t}\}_{t\in \mathbb{R}}$ at zero).

Further we note that, since $\left|e^{i t f(n)}-1\right| \leq\sqrt{2}|t||f(n)|,$ $n\in\mathbb{N},$ then
\begin{equation}\label{SS}
\left(\widehat{\Upsilon}_k^N(t)\right)^2= \sum\limits_{n=N+1}^{\infty} \left|\frac{e^{i t f(n)}-1}{t}\right|^2 \left|\Delta^{k} c_{n} \right|^2 \leq
2 \sum\limits_{n=N+1}^{\infty} \left|f(n)\Delta^{k} c_{n} \right|^2.
\end{equation}
Let's show that $\left\{ f(n)\Delta^{k} c_{n}\right\}_{n=1}^{\infty}\in \ell_2.$ Indeed, since by Leibnitz theorem~(\ref{Leibnitz})
$$\Delta^k (c_n f(n)) = \sum\limits_{j=0}^{k} C_{k}^j \Delta^{k-j} c_{n-j} \Delta^j f(n)=\sum\limits_{j=0}^{k}  C_{k}^j v_j(n),$$
and $\{f(n) \cdot c_{n}\}_{n=1}^{\infty}\in \ell_2(\Delta^k)$ from the condition~(\ref{Domain_k}), then $v= \sum\limits_{j=0}^{k} C_{k}^j v_j  \in \ell_2,$ where  $v_j=\{v_j(n)= \Delta^{k-j} c_{n-j} \Delta^j f(n)\}_{n=1}^{\infty}$ for $0\leq j\leq k.$ Using that $\left\{f(n)\right\}_{n=1}^{\infty}\in\mathcal{S}_k$ and applying the discrete Hardy inequality~(\ref{Hardy p=2}), for every $j:$ $1\leq j\leq k$ we will have
\begin{align*}
\|v_j\|_{\ell_2}^2 &= \sum\limits_{n=1}^{\infty} \left|\Delta^{k-j} c_{n-j} \Delta^j f(n) \right|^2=\sum\limits_{n=1}^{\infty} \frac{\left| n^{j}\Delta^j f(n) \right|^2 \left|\Delta^{k-j} c_{n-j} \right|^2}{n^{2j}}\\
&\leq C(j)^2\sum\limits_{n=1}^{\infty} \frac{\left|\Delta^{k-j} c_{n-j} \right|^2}{n^{2j}}= C(j)^2 \Omega_j\leq 4^j  M C(j)^2 \|x\|_k^{2}.
\end{align*}
(see estimate for $\Omega_j$ in the proof of the boundedness of $\{e^{A_k t}\}_{t\in \mathbb{R}}$).
Therefore
$v_0=\{f(n) \Delta^k c_n\}_{n=1}^{\infty}\in \ell_2$ and the remainder of the series in the right hand side of~(\ref{SS}) tends to zero when $N\rightarrow \infty$.

Since $\mathcal{Q}_k(|t|)$ is a polynomial of argument $|t|$, then there exists a constant $K\geq 0$ such that for all $t: \: |t|\leq 1$  we will have the following:
$$\mathcal{Q}_k(|t|)\leq K.$$
So, since $\{f(n) \Delta^k c_n\}_{n=1}^{\infty}\in \ell_2$ и $\{\Delta^k c_n\}_{n=1}^{\infty}\in \ell_2$, then for all $t: \: |t|\leq 1$ we will obtain that
\begin{align*}
  &\sum\limits_{n=N+1}^{\infty} \left|\Delta^k\left(\frac{e^{i t f(n)}-1}{t}c_{n}\right)\right|^2 \leq \left(\widehat{\Upsilon}_k^N(t) + K \left( \sum\limits_{n=N+1}^{\infty} \left|\Delta^k c_n \right|^2\right)^{\frac{1}{2}}\right)^2\\
  & \leq \left(\sqrt{2} \left(\sum\limits_{n=N+1}^{\infty} \left|f(n)\Delta^{k} c_{n} \right|^2\right)^{\frac{1}{2}} + K \left( \sum\limits_{n=N+1}^{\infty} \left|\Delta^k c_n \right|^2\right)^{\frac{1}{2}}\right)^2 \rightarrow 0
\end{align*}
when $N\rightarrow \infty$.

Further on, take any $\varepsilon>0$ and fix it.
Let $|t|\leq1.$ Since the condition~(\ref{first sum}) holds, taking into account the last relation, there exists $N\in \mathbb{N}$, such that
$$\sum\limits_{n=N+1}^{\infty} \left|\Delta^k\left(i f(n)c_{n}\right)\right|^2 \leq \varepsilon,\qquad \sum\limits_{n=N+1}^{\infty}  \left|\Delta^k\left(\frac{e^{i t f(n)}-1}{t}c_{n}\right)\right|^2 \leq \varepsilon.$$
It follows that
$\left\| \frac{e^{A_k t}x - x}{t} -A_k x\right\|_k \leq \frac{1}{\sqrt{m}} \Biggl(\sum\limits_{n=1}^{N} \left|\Delta^k\left(\left(\frac{e^{i t f(n)}-1}{t} - i f(n)\right)c_{n}\right)\right|^2+ 4  \varepsilon \Biggr)^{\frac{1}{2}}.
$
Since $\left\| \Delta^k\right\|\leq 2^k$, where $\Delta^k$ is considered as an operator from $\ell_2$ to $\ell_2,$
then
$$\left\| \frac{e^{A_k t}x - x}{t} -A_k x\right\|_k \leq \frac{1}{\sqrt{m}} \Biggl(2^{2k}\sum\limits_{n=1}^{N} \left|\left(\frac{e^{i t f(n)}-1}{t} - i f(n)\right)c_{n}\right|^2+ 4  \varepsilon \Biggr)^{\frac{1}{2}}.$$

Since for each $n\in[1,N]$, $\left|\left(\frac{e^{i t f(n)}-1}{t} - i f(n)\right)c_{n}\right|^2\rightarrow 0$, as $t\rightarrow 0$,
then $$2^{2k} \sum\limits_{n=1}^{N} \left|\left(\frac{e^{i t f(n)}-1}{t} - i f(n)\right)c_{n}\right|^2\rightarrow 0,$$ as $t\rightarrow 0$.
Therefore for earlier chosen $\varepsilon>0$ there exists $\delta(\varepsilon)>0$, such that for all $t:$ $|t|\leq \delta$ we have
$$2^{2k} \sum\limits_{n=1}^{N} \left|\left(\frac{e^{i t f(n)}-1}{t} - i f(n)\right)c_{n}\right|^2 \leq \varepsilon.$$
The last inequality also holds for $t:\:|t|\leq \min\{1,\delta\}.$
So
$$\left\| \frac{e^{A_k t}x - x}{t} -A_k x\right\|_k \leq \frac{1}{\sqrt{m}} \left(\varepsilon+ 4  \varepsilon \right)^{\frac{1}{2}}= \frac{\sqrt{5\varepsilon}}{\sqrt{m}}.$$
It means that $x\in E_k$ and $\left\| \frac{e^{A_k t}x - x}{t} -A_k x\right\|_k \rightarrow 0,$ as $t\rightarrow 0$.

Thus $A_k$ is a generator of constructed $C_0$-group $\{e^{A_k t}\}_{t\in \mathbb{R}}$ and the theorem is completely proved.
\hfill$\Box$
\end{pf}

We note that Theorem \ref{not-basis-ln} is a special case of Theorem \ref{not-basis_fk} when $k=1$ and $f(n)=\ln n,$ $n\in\mathbb{N}.$

Now we pass to the question on the asymptotic behaviour of constructed $C_0$-group $\{e^{A_k t}\}_{t\in \mathbb{R}}$.
\begin{prop}\label{grows as fk}
Let $k\in\mathbb{N}$ and suppose that $\left\{e^{A_k t} \right\}_{t\in \mathbb{R}}$ is a $C_0$-group from Theorem~\ref{not-basis_fk}. Then the following assertions hold true:
\begin{enumerate}
\item $\left\|e^{A_k t} \right\| \rightarrow \infty,$ as $t\rightarrow \pm\infty.$
\item There exists a polynomial $\mathfrak{p}_k$ with positive coefficients, with degree $\deg \mathfrak{p}_k =k,$ such that for every $t\in\mathbb{R}$ the following estimate is valid:
$$\left\|e^{A_k t}\right\| \leq \mathfrak{p}_k (|t|).$$
\end{enumerate}
\end{prop}
\begin{pf}
1. Recall that $A_k$ has eigenvalues $\{i f(n)\}_{n=1}^{\infty}\subset i\mathbb{R}$ and corresponding eigenvectors $\{e_n\}_{n=1}^{\infty}$ are dense in
$H_k\left(\{e_n\}\right)$.
By Theorem~\ref{not-basis_fk} operator $A_k$ generates the $C_0$-group $\left\{e^{A_k t} \right\}_{t\in \mathbb{R}}$.
Assume the opposite, namely that $\left\{e^{A_k t} \right\}_{t\in \mathbb{R}}$ is bounded $C_0$-group.
Then by theorem 2 from the paper of A.I.~Miloslavskii~\cite{Miloslavskii} the sequence $\{e_n\}_{n=1}^{\infty}$ forms a Riesz basis of $H_k\left(\{e_n\}\right)$.
Thus we arrive at a contradiction, since by
Proposition~\ref{prop} the sequence $\{e_n\}_{n=1}^{\infty}$ does not form a Schauder basis.

2. The existence of the polynomial $\mathfrak{p}_k$  with the desired properties for the case $k=1$ follows from the  proof of Theorem~\ref{not-basis-ln}, see estimate~(\ref{bounded ubove}), and for the case when $k>1$ it follows from the proof of Theorem~\ref{not-basis_fk}, see estimate~(\ref{bounded ubove k}).

\end{pf}

In 1967 V.E. Katsnel'son \cite{Katsnel'son} proved the following theorem.
\begin{thm}\cite{Katsnel'son}\label{Katsnel'son thm}
Let $\{\lambda_n\}_{n=1}^{\infty}$ be a sequence of distinct points in the upper half-plane $\{z\in \mathbb{C}:\: Im(z)>0\}$ and assume that
\begin{equation}\label{inf}
    \inf\limits_{1\leq j<\infty} \prod\limits_{k=1;k\neq j}^{\infty} \left|\frac{\lambda_j - \lambda_k}{\lambda_j - \overline{\lambda}_k}\right|=0.
\end{equation}
Then there exists a linear operator $A: H \supset D(A) \mapsto H$ such that:
\begin{enumerate}
\item $Im \langle Ax,x \rangle\geq 0,$ $x\in D(A).$
\item The eigenvalues of $A$ are $\{\lambda_n\}_{n=1}^{\infty}$ and $\sigma(A)=\overline{\{\lambda_n\}_{n=1}^{\infty}}$.
\item The system of eigenvectors $\{v_n\}_{n=1}^{\infty}$ of $A$ is dense in $H$ but not uniformly minimal.

Furthermore, if $\{Im(\lambda_n)\}_{n=1}^{\infty}$ is, in addition, a bounded sequence, then there exists a linear operator $A: H \supset D(A) \mapsto H$ satisfying 1-3 and
\item $A=A_{\Re} + i A_{\Im},$

where $A_{\Re}$ is selfadjoint operator and $A_{\Im}$ is bounded positive operator.
\end{enumerate}
\end{thm}
\begin{rmk}\label{Katsnel'son}\hspace{1mm}
\begin{itemize}
\item Note that in any horizontal strip $\{z\in \mathbb{C}:\: 0<Im(z)<\alpha\}$ the condition (\ref{inf}) turns into $\inf\limits_{j\neq k} |\lambda_j-\lambda_k|=0$, see \cite{Garnett}. It follows that for any sequence $\{\mu_n\}_{n=1}^{\infty}$ of distinct points in a vertical strip $\{z\in \mathbb{C}:\: -\alpha<Re(z)<0\}$ not satisfying (\ref{2}) there exists, by Theorem~\ref{Katsnel'son thm}, a dissipative operator $V=i A$ with eigenvalues $\{\mu_n\}_{n=1}^{\infty}$ and eigenvectors $\{v_n\}_{n=1}^{\infty}$, that are dense but not uniformly minimal in $H$. Consequently, by the Lumer-Phillips theorem, $V$ is an infinitesimal generator of contractive $C_0$-semigroup $\left\{T(t) \right\}_{t\geq 0}$.
    In this context we note that the $C_0$-semigroup $\left\{e^{A_k t} \right\}_{t\geq 0}$ constructed in Theorem~\ref{not-basis_fk} (the restriction of the $C_0$-group $\left\{e^{A_k t} \right\}_{t\in \mathbb{R}}$) is not contractive and grows when $t\rightarrow \infty,$ see Proposition~\ref{grows as fk}.
\item Since $\{Re(\mu_n)\}_{n=1}^{\infty}$ is a bounded sequence, the contractive semigroup $\left\{T(t) \right\}_{t\geq 0}$ expands to the $C_0$-group $\left\{T(t) \right\}_{t\in \mathbb{R}}$.
\item Theorem~\ref{not-basis_fk} can be easily generalized to the case when $\sigma(A_k)$ belongs to arbitrary vertical line. But Theorem~\ref{not-basis_fk} is focused on a critical for Theorem~\ref{Katsnel'son thm} case when $\sigma(A_k)\subset i\mathbb{R}.$
\item The construction of unbounded generator of a $C_0$-group on $H$ with non-bounded non-Riesz basis family of eigenvectors is trivial.
However, the existence of unbounded generator of the $C_0$-group on $H$ with bounded non-Riesz basis family of eigenvectors is unknown.
\end{itemize}
\end{rmk}

The proof of Theorem \ref{not-basis_fk} leads to the following.
\begin{corollary}\label{Cor1}
Let $\omega_{0,k}$ is the growth bound of the $C_0$-group $e^{A_k t}$ constructed in Theorem~\ref{not-basis_fk} and $s(A_k)=\sup\{Re \lambda: \: \lambda\in\sigma(A_k)\}$ is a spectral bound of its generator $A_k$. Then, for any $k\in\mathbb{N},$
$$\omega_{0,k}=s(A_k)=0.$$
\end{corollary}

Combining Theorem~\ref{not-basis_fk} and Corollary~\ref{Cor1} with the spectral theorem of K.~Boyadzhiev and R.~DeLaubenfels~\cite{Boyadzhiev} we obtain the following result.
\begin{prop}\label{p1}
For each $k\in\mathbb{N}$ and arbitrary $\alpha>0$ the operator $A_k$ constructed in Theorem~\ref{not-basis_fk} has a bounded $\mathcal{H}^{\infty}$-calculus on a strip $H_\alpha=\{z\in \mathbb{C}:\: |Re(z)|<\alpha\}$.
\end{prop}
It is interesting to compare this result with the construction of an operator $A$ on $H$ without a bounded $\mathcal{H}^{\infty}$-calculus from~\cite{Haase2}(Section 5.5). This construction is based on the concept of non-Riesz basis. About $\mathcal{H}^{\infty}$-calculus see~\cite{Boyadzhiev,Haase1,Haase2,Zwart}, for example.

Within the context of Theorem \ref{not-basis_fk} we finally note the following.
\begin{rmk}\label{Remark}\hspace{1mm}
\begin{itemize}
\item For each $k\in\mathbb{N}$ the operator $A_k$ from Theorem \ref{not-basis_fk} is an unbounded linear operator with $\overline{D(A_k)}=H_k\left(\{e_n\}\right)$
and $A_k$ is closed on $D(A_k)$.
\item $A_k$ has pure point spectrum $\{i f(n)\}_{n=1}^{\infty}$, which does not satisfy the condition~(\ref{2}), and, moreover, cannot be decomposed into $K$ sets, with every set having a uniform gap, since $\{f(n)\}_{n=1}^{\infty}\in\mathcal{S}_k$.
\item If we omit the condition (\ref{2}), then the converse, in some sense, statement to the Theorem~\ref{Remarkable result} holds. More precisely, suppose that $f: [1,+\infty)\mapsto \mathbb{R}$ is any real function such that $\{f(n)\}_{n=1}^{\infty}$ does not satisfy the condition (\ref{2}). If we define $A: H \supset D(A) \mapsto H$ as $A \sum\limits_{n=1}^{\infty} \alpha_n e_n = \sum\limits_{n=1}^{\infty} i f(n) \cdot \alpha_n e_n,$ with domain
    $$D(A)=\left\{x= \sum\limits_{n=1}^{\infty} \alpha_n e_n \in H:\: \left\{f(n) \cdot \alpha_n\right\}_{n=1}^{\infty} \in \ell_2 \right\},$$ then $A$ is a Riesz-spectral operator and it generates a $C_0$-group on $H$. For details see~\cite{Curtain}.
\end{itemize}
\end{rmk}

\section{The construction of infinitesimal operators with non-basis family of eigenvectors on certain Banach spaces}
The construction of infinitesimal operators with complete minimal non-basis family of eigenvectors on Banach spaces
is similar to the construction of infinitesimal operators with complete minimal non-basis family of eigenvectors on Hilbert spaces. Namely, we have the following theorem, analogous to Theorem~\ref{not-basis_fk}.
\begin{thm}\label{not-basis-ell_p_fk}
Assume that $\{e_n\}_{n=1}^{\infty}$ is a symmetric basis in $\ell_p,$ $p>1$, and $k\in\mathbb{N}.$ Then $\{e_n\}_{n=1}^{\infty}$ does not form a basis of $\ell_{p,k}\left(\{e_n\}\right)$
and the operator $A_k:\ell_{p,k}\left(\{e_n\}\right) \supset D(A_k) \mapsto \ell_{p,k}\left(\{e_n\}\right),$ defined by
$$A_k x=A_k (\mathfrak{f})\sum\limits_{n=1}^{\infty} c_{n} e_n= (\mathfrak{f})\sum\limits_{n=1}^{\infty} i f(n) \cdot c_{n} e_n,$$
where $\left\{f(n)\right\}_{n=1}^{\infty}\in\mathcal{S}_k$, with domain
\begin{equation}\label{Domain_k_ell_p}
    D(A_k)=\left\{x= (\mathfrak{f})\sum\limits_{n=1}^{\infty} c_{n} e_n \in \ell_{p,k}\left(\{e_n\}\right):\:  \{f(n) \cdot c_{n}\}_{n=1}^{\infty}\in \ell_p(\Delta^k)\right\},
\end{equation}
generates the $C_0$-group on $\ell_{p,k}\left(\{e_n\}\right)$, which acts for every $t\in\mathbb{R}$ by the formula
$$e^{A_k t} x=e^{A_k t}(\mathfrak{f})\sum\limits_{n=1}^{\infty} c_n e_n =(\mathfrak{f})\sum\limits_{n=1}^{\infty} e^{i t f(n)} c_{n} e_n.
$$
\end{thm}
\begin{pf} The proof of Theorem~\ref{not-basis-ell_p_fk} is based on the combination of Proposition~\ref{props} and Proposition~\ref{sym} with Hardy inequality~(\ref{1}) for the case $p>1$ and may be performed similarly to the proof of Theorem~\ref{not-basis_fk}.
\hfill$\Box$
\end{pf}

Concerning the asymptotic behaviour of $C_0$-group $\left\{e^{A_k t} \right\}_{t\in \mathbb{R}}$ constructed in the Theorem~\ref{not-basis-ell_p_fk} we have the following statement, analogous to the point 2 of Proposition~\ref{grows as fk}.
\begin{prop}\label{grows as ell p fk}
Let $k\in\mathbb{N}$ and suppose that $\left\{e^{A_k t} \right\}_{t\in \mathbb{R}}$ is a $C_0$-group from Theorem~\ref{not-basis-ell_p_fk}. Then there exists a polynomial $\mathfrak{p}_k$ with positive coefficients, with degree $\deg \mathfrak{p}_k =k,$ such that for every $t\in\mathbb{R}$ the following estimate is valid:
$$\left\|e^{A_k t}\right\| \leq \mathfrak{p}_k (|t|).$$
\end{prop}

So we see that our class of $C_0$-groups belong to the class of polynomially bounded $C_0$-groups.

Define the operator $B: \ell_{p,1}\left(\{e_n\}\right) \supset D(B) \mapsto \ell_{p,1}\left(\{e_n\}\right)$ as
\begin{equation}\label{B}
  B x=B (\mathfrak{f})\sum\limits_{n=1}^{\infty} c_n e_n= (\mathfrak{f})\sum\limits_{n=1}^{\infty} \lambda_n c_n e_n,
\end{equation}
with domain
\begin{equation}\label{Domain_B}
    D(B)=\left\{x= (\mathfrak{f})\sum\limits_{n=1}^{\infty} c_n e_n \in \ell_{p,1}\left(\{e_n\}\right):\:\{\lambda_n c_n\}_{n=1}^{\infty}\in \ell_p(\Delta)\right\},
\end{equation}
In particular case when $k=1$ and $f(x)=\ln x$ we have from Theorem~\ref{not-basis-ell_p_fk} an immediate consequence, similar to Theorem~\ref{not-basis-ln}.
\begin{corollary}\label{not-basis-ell_p-ln}
Let $\{e_n\}_{n=1}^{\infty}$ be a symmetric basis of $\ell_p,$ $p>1$. Then $\{e_n\}_{n=1}^{\infty}$ does not form a basis of $\ell_{p,1}\left(\{e_n\}\right)$ and the operator $B$ defined by (\ref{B}) with domain (\ref{Domain_B}), where $\lambda_n=i\ln n$,
generates a $C_0$-group on $\ell_{p,1}\left(\{e_n\}\right)$.
\end{corollary}

As in the case of Proposition~\ref{sqrt} we note that even if we consider the spectrum $\{\lambda_n\}_{n=1}^{\infty}$ of $B$ defined by (\ref{B},\ref{Domain_B}) of the same geometric nature, i.e. satisfying~(\ref{geometric cond}), then $B$ not necessarily generates a $C_0$-semigroup on $\ell_{p,1}\left(\{e_n\}\right)$.
Using the same arguments as in the proof of Proposition~\ref{sqrt} we can obtain the following.

\begin{prop}\label{sqrt h p}
Let $p\geq 1$, and suppose that $\{\lambda_n\}_{n=1}^{\infty}\subset i\mathbb{R}$ satisfies~(\ref{geometric cond}) and, moreover, there exists $\alpha\in \left(0,\frac{1}{p}\right]:$ $\liminf\limits_{n\rightarrow \infty} n^{\alpha} |\lambda_{n-1}-\lambda_n|>0.$
Then the operator $A$, defined by~(\ref{B}), with domain~(\ref{Domain_B}), does not generate a $C_0$-semigroup on $\ell_{p,1}\left(\{e_n\}\right)$.
\end{prop}

Next, we deduce a consequence, similar to Corollary~\ref{Cor1}.
\begin{corollary}\label{Cor2}
Let $\omega_{0,k}$ is the growth bound of the $C_0$-group $e^{A_k t}$ constructed in Theorem~\ref{not-basis-ell_p_fk} and $s(A_k)=\sup\{Re \lambda: \: \lambda\in\sigma(A_k)\}$ is a spectral bound of its generator $A_k$. Then, for any $k\in\mathbb{N},$
$$\omega_{0,k}=s(A_k)=0.$$
\end{corollary}
We complete the section with the following observation.
\begin{rmk}\label{Remark_ell_p}\hspace{1mm}
\begin{itemize}
\item For each $k\in\mathbb{N}$ the operator $A_k$ from Theorem \ref{not-basis-ell_p_fk} is an unbounded linear operator with $\overline{D(A_k)}=\ell_{p,k}\left(\{e_n\}\right)$ and $A_k$ is closed on $D(A_k)$.
\item $A_k$ has pure point spectrum $\{i f(n)\}_{n=1}^{\infty}$, which does not satisfy the condition~(\ref{2}), and, moreover, cannot be decomposed into $K$ sets, with every set having a uniform gap, since $\{f(n)\}_{n=1}^{\infty}\in\mathcal{S}_k$.
\item Suppose that $f: [1,+\infty)\mapsto \mathbb{R}$ is any real function. If we define $A: \ell_p \supset D(A) \mapsto \ell_p$ as $A \sum\limits_{n=1}^{\infty} \alpha_n e_n = \sum\limits_{n=1}^{\infty} i f(n) \cdot \alpha_n e_n,$ with domain $$D(A)=\left\{x= \sum\limits_{n=1}^{\infty} \alpha_n e_n \in \ell_p:\: \left\{f(n)\cdot \alpha_n\right\}_{n=1}^{\infty} \in \ell_p \right\},$$ then it can be shown that $A$ generates a bounded $C_0$-group on $\ell_p$.
\end{itemize}
\end{rmk}
\section{Concluding remarks}
The results of the present paper allow us to say the following. A Theorem~\ref{Remarkable result} cannot be improved. Moreover, it is impossible to obtain any analogue of Theorem \ref{Remarkable result} concerning non-basis family of eigenvectors by means of omitting or weakening of the condition (\ref{2}). On the other hand, it is interesting to obtain some analogues of Theorem \ref{Remarkable result} in Banach spaces with certain classes of bases, e.g., symmetric bases, unconditional bases.

Theorem \ref{not-basis-ell_p_fk} and Remark \ref{Remark_ell_p} allow us to say that symmetric bases in $\ell_p$ spaces behave like Riesz bases in $H$. Consequently, we arrive at the idea of possible generalization of the Theorem~\ref{Remarkable result} to the case of operators with symmetric basis family of eigenvectors, that generate $C_0$-groups on the spaces $\ell_p,\:p\geq 1$, and propose the following.
\begin{conjecture}\label{conjecture}
Let $A$ be the generator of the $C_0$-group on the space $\ell_p,\:p\geq 1,$ with eigenvalues $\{\lambda_n\}_{n=1}^{\infty}$ (counting with multiplicity) and the corresponding (normalized) eigenvectors $\{e_n\}_{n=1}^{\infty}$.

If $\overline{Lin}\{e_n\}_{n=1}^{\infty}=\ell_p$ and the point spectrum $\{\lambda_n\}_{n=1}^{\infty}$ satisfies (\ref{2}), then $\{e_n\}_{n=1}^{\infty}$ forms a symmetric basis of $\ell_p$.
\end{conjecture}

Finally, the answer to the following important question does not yet exist. Is it possible to construct the unbounded generator of a $C_0$-group with bounded non-Riesz basis family of eigenvectors?
We note that the construction of unbounded generator of a $C_0$-semigroup with bounded non-Riesz basis family of eigenvectors is quite simple, see, e.g.,~\cite{Haase2}(Section 5.5).

\bibliography{mybibfile}

\end{document}